\documentclass[a4paper,10pt]{amsart}

\usepackage[top=2.5cm, bottom=2.5cm, left=2.5cm, right=2.5cm]{geometry}

\usepackage{amssymb,bm}
\usepackage[dvipdfmx]{graphicx}
\usepackage{xcolor}

\newtheorem{thm}{Theorem}[section]

\newtheorem{rem}[thm]{Remark}

\numberwithin{equation}{section}

\allowdisplaybreaks

\newcommand{\mcd}{\mathcal{D}}
\newcommand{\mcf}{\mathcal{F}}
\newcommand{\mci}{\mathcal{I}}

\newcommand{\mcl}{\mathcal{L}}

\newcommand{\mfL}{\mathfrak{L}}

\newcommand{\mbbh}{\mathbb{H}}

\newcommand{\mbbn}{\mathbb{N}}

\newcommand{\mbbr}{\mathbb{R}}

\newcommand{\mbbrp}{\mathbb{R}_{+}}
\newcommand{\mbbzp}{\mathbb{Z}_{+}}

 \newcommand{\lam}{\lambda} \newcommand{\ep}{\epsilon} 
\newcommand{\vp}{\varphi} \newcommand{\del}{\delta} \newcommand{\sig}{\sigma}
\newcommand{\D}{\Delta}  
\newcommand{\Lam}{\Lambda} \newcommand{\Gam}{\Gamma}

\newcommand{\p}{\partial}
\newcommand{\ra}{\rangle} \newcommand{\la}{\langle}
\newcommand{\wc}{\Rightarrow} 
\newcommand{\cil}{\xrightarrow{\mcl}} 
\newcommand{\scl}{\xrightarrow{\mcl_{s}}} 
\newcommand{\cip}{\xrightarrow{p}} 

\newcommand{\argmin}{\mathop{\rm argmin}}

\def\ds#1{\displaystyle{#1}}

\newcommand{\hmcheck}{}  

\def\nn{\nonumber}

\def\v2#1{\textcolor{black}{#1}}
\newcommand{\pr}{P} \newcommand{\E}{E}

\def\diag{{\rm diag}}
\def\sumj{\sum_{j=1}^{n}}

\def\lp{L\'evy process}

\def\cadlag{c\`adl\`ag}

\newcommand{\tes}{\hat{\theta}_{n}}
\newcommand{\tesnn}{\hat{\theta}}
\newcommand{\mesnn}{\hat{\mu}}
\newcommand{\lesnn}{\hat{\lambda}}
\newcommand{\besnn}{\hat{\beta}}
\newcommand{\sesnn}{\hat{\sigma}}

\newcommand{\tz}{\theta_0}

\newcommand{\sgn}{\mathrm{sgn}}

\title[Optimal stable OU regression]
{Optimal stable Ornstein-Uhlenbeck regression}

\author[H. Masuda]{Hiroki Masuda}
\address{Faculty of Mathematics, Kyushu University, 744 Motooka Nishi-ku Fukuoka 819-0395, Japan
\and 
Graduate School of Mathematical Sciences, The University of Tokyo, 3-8-1 Komaba Meguro-ku Tokyo 153-8914, Japan.
}
\email{hmasuda@ms.u-tokyo.ac.jp}

\date{}

\begin{document}
\setlength{\baselineskip}{4.5mm}

\maketitle

\begin{abstract}
We prove asymptotically efficient inference results concerning an Ornstein-Uhlenbeck regression model driven by a non-Gaussian stable {\lp}, where the output process is observed at high frequency over a fixed period. The local asymptotics of non-ergodic type for the likelihood function is presented, followed by a way to construct an asymptotically efficient estimator through a suboptimal yet very simple preliminary estimator.
\end{abstract}

%
%

\section{Introduction}

\subsection{Objective and background}

Given an underlying filtered probability space $(\Omega,\mcf,(\mcf_t)_{t\in[0,T]},\pr)$, we consider the following Ornstein-Uhlenbeck (OU) regression model
\begin{equation}
Y_t = Y_0 +\int_0^t (\mu\cdot X_s-\lam Y_s)ds + \sig J_t, \qquad t\in[0,T],
\label{soureg:model}
\end{equation}
where $J$ is the symmetric $\beta$-stable ({\cadlag}) {\lp} characterized by
\begin{equation}
\E[e^{iu J_t}] = \exp(-t|u|^{\beta}),\qquad t\ge 0,~u\in\mbbr,
\nn
\end{equation}
and is independent of the initial variable $Y_0$,
and where $X=(X_{t})_{t\in[0,T]}$ is an $\mbbr^{q}$-valued non-random {\cadlag} function such that
\begin{equation}
\lam_{\min}\left(\int_0^T X_t^{\otimes 2}dt \right) >0,
\label{A_X.pd}
\end{equation}
\v2{with $\lam_{\min}(A)$ denoting the minimum eigenvalue of a square matrix $A$.}
Throughout, the terminal sampling time $T>0$ is a fixed constant.
Let
\begin{equation}
\theta := (\lam,\mu,\beta,\sig) \in \Theta, 
\nonumber
\end{equation}
where $\Theta\subset\mbbr^{p}$ ($p:=q+3$) is a bounded convex domain such that its closure $\overline{\Theta}\subset \mbbr\times\mbbr^q \times(0,2)\times(0,\infty)$.
%
The primary objective of this paper is the asymptotically efficient estimation of $\theta$, when available data is $(X_t)_{t\in[0,T]}$ and $(Y_{t_{j}})_{j=1}^n$, where $t_j=t_j^n:=jh$ with $h=h_n:=T/n$;
later on, we will consider cases where we observe $(X_{t_{j}})_{j=1}^n$ instead of the full continuous-time record.
We will denote the true value of $\theta$ by $\tz=(\lam_0,\mu_0,\beta_0,\sig_0)\in\Theta$.

Analysis of the (time-homogeneous) OU process driven by a stable {\lp} goes back to \cite{Doo42}, where Doob treated the model in a genuinely analytic manner without It\^o's formula, which has not got published as yet at that time.
Nowadays, the OU models have been used in a wide variety of applications, such as electric consumption modeling \cite{PKAS11}, \cite{BorSch17}, and \cite{VAKB19}, 
ecology \cite{JhwMar14}, and protein dynamics modeling \cite{ChaSch12}, to mention a few.

The model \eqref{soureg:model} may be seen as a continuous-time counterpart of the simple first-order ARX (autoregressive exogenous) model.
Nevertheless, any proper form of the efficient estimation result has been missing in the literature, probably due to the lack of background theory that can deal with the estimation of all the parameters involved under the bounded-domain infill asymptotics.
Let us note that, when $J$ is a Wiener process ($\beta=2$), the drift parameters are consistently estimable only when the terminal sampling time tends to infinity, and the associated statistical experiments are known to possess essentially different properties according to the sign of $\lam$. That is to say, the model is:
locally asymptotically normal for $\lam>0$ (ergodic case);
locally asymptotically Brownian functional for $\lam=0$ (unit-root case);
locally asymptotically mixed-normal (LAMN) for $\lam<0$ (non-ergodic (explosive) case).
Turning back to the stable driven case, we should note that the least-squares type estimator would not work unless $T_n\to\infty$, as is expected from \cite{HuLon09} and \cite{ZhaZha13}; there, the authors proved that (when $\beta$ is known) the rate of convergence when $\lam>0$ equals $(T_n/\log n)^{1/\beta}$ and the asymptotic distribution is given by a ratio of two independent stable distributions.

\subsection{Contributions in brief}

First, in Section \ref{sec:likelihood.asymptotics}, we will show that the model is locally asymptotically mixed normal (LAMN) at $\tz\in\Theta$, and also that the likelihood equation has a root that is asymptotically efficient in the classical sense of Haj\'{e}k-Le Cam-Jeganathan.
The asymptotic results presented here are uniformly valid in a single manner over any compact subset of the parameter space $\Theta$.
In particular, the sign of the autoregressive parameter $\lam_0$ does not matter, revealing that (i) the results can be described in a unified manner regardless of whether the model is ergodic or not, and also that (ii) the conventional unit-root problem (see \cite{SamKni09} and the references therein) is not relevant here at all;
this is in sharp contrast to the case of ARX time-series models and the Gaussian OU models.
Besides, in Section \ref{sec:efficient.estimator}, we will provide a way to provide an asymptotically efficient estimator through a suboptimal yet very simple preliminary estimator, which enables us to bypass not only computationally demanding numerical optimization of the likelihood function involving the $\beta$-stable density, but also possible multiple-root problem \cite[Section 7.3]{Leh99}.

\section{Local likelihood asymptotic}
\label{sec:likelihood.asymptotics}

\subsection{Preliminaries and result}

Let $\pr_\theta$ denote the image measure of $(J,Y)$ associated with the value $\theta\in\Theta$.
We are going to show the non-trivial stochastic expansion of the log-likelihood ratio of $\pr_{\theta+\vp_{n}(\theta)v_n,n}^Y$ with respect to $\pr_{\theta,n}^Y$ for appropriate norming matrix $\vp_n(\theta)$ introduced later and bounded sequence $(v_n)\subset\mbbr^p$, where $\pr_{\theta,n}^Y$ stands for the restriction of $\pr_\theta$ to $\sig(Y_{t_j}:\,j\le n)$.
The distribution $\mcl(Y_0)$ may vary according as $\theta$, we will assume that for any $\ep>0$ there exists an $M>0$ such that $\sup_{\theta\in\overline{\Theta}} \pr_\theta[|Y_0|\ge M]<\ep$.

Let $\phi_\beta$ denote the $\beta$-stable density of $J_1$: $\pr_\theta[J_1\in dy]=\phi_\beta(y)dy$.
It is known that $\phi_{\beta}(y)>0$ for each $y\in\mbbr$, that $\phi_{\beta}$ is smooth in $(y,\beta)\in\mbbr\times(0,2)$, and that for each $k,l\in \mbbzp$,
\begin{equation}
\limsup_{|y|\to\infty}
\frac{|y|^{\beta+1+k}}{\log^{l}(1+|y|)} \big|\p^{k}\p_{\beta}^{l}\phi_{\beta}(y)\big| < \infty.
\label{hm:dens.ab}
\end{equation}
See \cite{DuM73} for details.
Here we wrote $\p^{k}\p_{\beta}^{l}\phi_{\beta}(y):=(\p^{k}/\p y^k)(\p^l/\p\beta^{l})\phi_{\beta}(y)$; analogous notation for the partial derivatives will be used in the sequel.

To proceed, we need to introduce further notation.
Any asymptotics will be taken for $n\to\infty$ unless otherwise mentioned.
We denote by $\to_{u}$ the uniform convergence of non-random quantities concerning $\theta$ over $\overline{\Theta}$.
We write $C$ for positive universal constant which may vary at each appearance, and $a_{n}\lesssim b_{n}$ when $a_{n}\le C b_{n}$ for every $n$ large enough.
Given positive functions $a_{n}(\theta)$ and $b_{n}(\theta)$, we write $b_n(\theta)=o_u(a_n)$ and $b_n(\theta)=O_u(a_n)$ if $a_n^{-1}b_n(\theta)\to_u 0$ and $\sup_\theta |a_n^{-1}b_n(\theta)| =O(1)$, respectively. The symbol $a_{n}(\theta) \lesssim_{u} b_{n}(\theta)$ means that $\sup_{\theta}|a_{n}(\theta)/b_{n}(\theta)| \lesssim 1$. We write $\int_j$ instead of $\int_{t_{j-1}}^{t_j}$.

By integrating by parts applied to the process $t\mapsto e^{\lam t}Y_t$, we obtain the explicit {\cadlag} solution process: under $\pr_\theta$,
\begin{equation}
Y_t  = e^{-\lam (t-s)}Y_s + \mu\cdot\int_s^t e^{-\lam(t -u)} X_{u} du + \sig \int_s^t e^{-\lam(t -u)}dJ_u,\qquad t>s.
\label{Y_sol}
\end{equation}
For $x,\lam\in\mbbr$, we write
\begin{align}
\eta(x)=\frac{1}{x}(1-e^{-x}), \qquad 
\zeta_{j}(\lam) = \frac{1}{h}\int_{j} e^{-\lam(t_{j}-s)} X_s ds.
\nonumber
\end{align}
The basic property of the L\'{e}vy integral 
\v2{
and the fact that $\log \E_\theta[e^{iu J_1}]=-|u|^\beta$ give
}
\begin{align}
\log \E_\theta\left[\exp\left(iu\,\sig \int_{j} e^{-\lam(t_j -s)}dJ_s \right)\right]
&=
\v2{
\int_j \log \E_\theta\left[\exp\left( iu e^{-\lam(t_j -s)} \sig J_1\right)\right] ds
}
\nn\\
&= - |\sig u|^\beta \int_{j} e^{-\lam\beta(t_j -s)}ds 
\nn\\
&= - \big|\sig h^{1/\beta} \eta(\lam\beta h)^{1/\beta} u \big|^\beta.
\nonumber
\end{align}
Hence
\begin{equation}
\ep_{j}(\theta) := \frac{Y_{t_{j}} - e^{-\lam h}Y_{t_{j-1}} - \mu \cdot \zeta_{j}(\lam)h}{\sig h^{1/\beta}\eta(\lam\beta h)^{1/\beta}}
~\stackrel{\pr_{\theta}}{\sim}~\text{i.i.d.}~\mcl(J_1).
\label{epj_def}
\end{equation}
Now, the exact log-likelihood function $\ell_n(\theta)=\ell_n\left(\theta;\,(X_t)_{t\in[0,T]},(Y_{t_j})_{j=0}^n\right)$ is given by
\begin{align}
\ell_n(\theta) 
&=\sumj \log\left(\frac{1}{\sig h^{1/\beta}\eta(\lam\beta h)^{1/\beta}}\phi_{\beta} \left( \ep_{j}(\theta) \right)\right) \nn\\
&=\sumj \left( -\log\sig +\frac1\beta \log(1/h) 
\v2{- \frac1\beta \log\eta(\lam\beta h)} + \log\phi_{\beta} \left( \ep_{j}(\theta) \right)\right).
\label{log-lf}
\end{align}

\medskip

We introduce the non-random $p\times p$-matrix
\begin{equation}
\vp_n=\vp_{n}(\theta) := \diag\left( \frac{1}{\sqrt{n}h^{1-1/\beta}}\,I_{1+q},~
\frac{1}{\sqrt{n}}
\begin{pmatrix}
\vp_{11,n}(\theta) & \vp_{12,n}(\theta) \\
\vp_{21,n}(\theta) & \vp_{22,n}(\theta) \\
\end{pmatrix}
\right),
\label{hm:vp_def}
\end{equation}
where the real entries $\vp_{kl,n}=\vp_{kl,n}(\theta)$ are assumed to be continuously differentiable in $\theta \in \Theta$ and to satisfy the following conditions for some finite values $\overline{\vp}_{kl}=\overline{\vp}_{kl}(\theta)$:
\begin{equation}
\left\{
\begin{array}{l}
\vp_{11,n}(\theta) \to_{u} \overline{\vp}_{11}(\theta), \\[1mm]
\vp_{12,n}(\theta) \to_{u} \overline{\vp}_{12}(\theta), \\[1mm]
s_{21,n}(\theta):=\beta^{-2}\log(1/h_n)\vp_{11,n}(\theta) + \sig^{-1}\vp_{21,n}(\theta) \to_{u} \overline{\vp}_{21}(\theta), \\[1mm]
s_{22,n}(\theta):=\beta^{-2}\log(1/h_n)\vp_{12,n}(\theta) + \sig^{-1}\vp_{22,n}(\theta) \to_{u} \overline{\vp}_{22}(\theta), \\[2mm]
\ds{\inf_{\theta}|\overline{\vp}_{11}(\theta)\overline{\vp}_{22}(\theta) - \overline{\vp}_{12}(\theta)\overline{\vp}_{21}(\theta)|>0,} \\[1mm]
\ds{\max_{(k,l)} \left|\p_{\theta} \vp_{kl,n}(\theta)\right| \lesssim_u \log^2(1/h)}.
\end{array}
\right.
\label{hm:vp-conditions}
\end{equation}
The matrix $\vp_n(\theta)$ will turn out to be the right norming with which $u \mapsto \ell_{n}\left(\theta+\vp_{n}(\theta)u\right) - \ell_{n}\left(\theta\right)$ under $\pr_\theta$ has an asymptotically quadratic structure in $\mbbr^{p}$;
see \cite{BroMas18} and \cite{CleGlo20} for the related previous studies.
Note that $\sqrt{n}h_n^{1-1/\beta}\to_{u}\infty$ and $|\vp_{21,n}(\theta)|\vee|\vp_{22,n}(\theta)| \lesssim \log(1/h)$.
By the same reasoning as in \cite[page 292]{BroMas18}, we have $\inf_{\theta}|\vp_{11,n}(\theta)\vp_{22,n}(\theta) - \vp_{12,n}(\theta)\vp_{21,n}(\theta)| \gtrsim 1$ and $|\vp_{n}(\theta)| \to_{u} 0$ under \eqref{hm:vp-conditions}.

Let
\begin{equation}
f_\beta(y):= \frac{\p_{\beta}\phi_{\beta}}{\phi_{\beta}}(y), \qquad g_\beta(y) := \frac{\p\phi_{\beta}}{\phi_{\beta}}(y),
\nonumber
\end{equation}
and define the block-diagonal random matrix
\begin{equation}
\mci(\theta)=\diag\big(\mci_{\lam,\mu}(\theta),\mci_{\beta,\sig}(\theta)\big),
\label{hm:FIm_def0}
\end{equation}
where, for a random variable $\ep \stackrel{\pr_{\theta}}{\sim} \phi_\beta(y)dy$ and by denoting by $A^\top$ the transpose of matrix $A$,
\begin{align}
\mci_{\lam,\mu}(\theta)
&:= \frac{1}{\sig^2}\E_\theta\left[ g_\beta(\ep)^2 \right]
\frac1T\int_0^T
\begin{pmatrix}
Y_t^2 & -Y_t X_t^{\top} \\
-Y_t X_t & X_t^{\otimes 2}
\end{pmatrix}
dt, \label{hm:FIm_def1}\\
\mci_{\beta,\sig}(\theta)
&:= \begin{pmatrix}
\overline{\vp}_{11} & \overline{\vp}_{12} \\
-\overline{\vp}_{21} & -\overline{\vp}_{22}
\end{pmatrix}^{\top}\!\!\!
\begin{pmatrix}
\E_\theta\left[f_\beta(\ep)^2\right] & \E_\theta\left[\ep f_\beta(\ep)g_\beta(\ep)\right] \\
\E_\theta\left[\ep f_\beta(\ep)g_\beta(\ep)\right] & \E_\theta\left[(1+\ep g_\beta(\ep))^{2}\right]
\end{pmatrix}\!\!
\begin{pmatrix}
\overline{\vp}_{11} & \overline{\vp}_{12} \\
-\overline{\vp}_{21} & -\overline{\vp}_{22}
\end{pmatrix}.
\label{hm:FIm_def2}
\end{align}
Note that $\mci(\theta)$ does depend on the choice of $\overline{\vp}(\theta)=\{\overline{\vp}_{kl}(\theta)\}$;
if $\overline{\vp}(\theta)$ is free from $(\lam,\mu)$, then so is $\mci(\theta)$.

Also, we note that $\mci(\theta) >0$ ($\pr_\theta$-a.s., $\theta\in\Theta$) under \eqref{A_X.pd}.
Indeed, it was verified in \cite[Theorem 1]{BroMas18} that $\mci_{\beta,\sig}(\theta)>0$ a.s.
To deduce that $\mci_{\lam,\mu}(\theta)>0$ a.s., we note that $\int_0^T Y^2_t dt>0$ a.s. and that, by Schwarz's inequality,
\begin{align}
& u^\top \left\{ \int_0^T X_t^{\otimes 2} dt - 
\left( \int_0^T Y_t X_tdt \right) \left( \int_0^T Y_t^2 dt \right)^{-1} \left( \int_0^T Y_t X_tdt \right)^\top\right\} u
\nn\\
&= \int_0^T (u\cdot X_t)^2 dt - \left( \int_0^T Y_t^2 dt \right)^{-1} \left( \int_0^T Y_t (u\cdot X_t)dt \right)^2 > 0
\nonumber
\end{align}
for every nonzero $u\in\mbbr^q$, since for any constant real $\xi$ we have $Y\ne (u\cdot X)\xi$ a.s. as functions on $[0,T]$.
Apply the identity $\det\begin{pmatrix} A & B^\top \\ B & C\end{pmatrix}=\det(A)\det (C-BA^{-1}B^\top)$ to conclude the $\pr_\theta$-a.s. positive definiteness of $\mci(\theta)$.

The normalized score function $\D_n(\tz)$ and the normalized observed information matrix $\mci_n(\tz)$ are given by
\begin{align}
\D_{n}(\theta) &:= \vp_{n}(\theta)^{\top}\p_{\theta}\ell_{n}(\theta), \nn\\
\mci_n(\theta) &:= -\vp_n(\theta)^{\top}\p_{\theta}^{2}\ell_n(\theta)\vp_n(\theta),
\nonumber
\end{align}
respectively.
Let $MN_{p,\theta}(0,\mci(\theta)^{-1})$ denote the covariance mixture of $p$-dimensional normal distribution, corresponding to the characteristic function $u\mapsto\E_\theta\left[\exp(-u^\top\mci(\theta)^{-1}u/2)\right]$.
Finally, we write $M[u]=\sum_i M_i u_i$ for a linear form $M=\{M_i\}$ and similarly $Q[u,u]=Q[u^{\otimes 2}]=\sum_{i,j}Q_{ij}u_i u_j$ for a quadratic form $Q=\{Q_{ij}\}$.
Now we are ready to state the main claim of this section.

\begin{thm}
\label{thm_local.LF.asymp}
The following statements hold for any $\theta\in\Theta$. 
\begin{enumerate}
\item 
For any bounded sequence $(v_n)\subset\mbbr^p$, it holds that
\begin{equation}
\ell_{n}\left(\theta+\vp_{n}(\theta)v_n\right) - \ell_{n}\left(\theta\right) 
= \D_n(\theta)[v_n] - \frac{1}{2} \mci_n(\theta) [v_n,v_n] + o_{\pr_{\theta}}(1),
\nn
\end{equation}
where we have the convergence in distribution under $\pr_\theta$: 
$\mcl\left( \D_{n}(\theta), \, \mci_n(\theta) |\pr_\theta\right) \wc \mcl\left( \mci(\theta)^{1/2}Z,\, \mci(\theta) \right)$,
where $Z\sim N_{p}(0,I)$ is independent of \v2{$\mci(\theta)$}, defined on an extended probability space.
\item There exists a local maximum point $\tes$ of $\ell_{n}(\theta)$ with $\pr_\theta$-probability tending to $1$ for which
\begin{equation}
\vp_{n}(\theta)^{-1}(\tes -\theta) = \mci_{n}(\theta)^{-1}\D_n(\theta) + o_{\pr_\theta}(1)
\wc MN_{p,\theta}\left(0,\, \mci(\theta)^{-1} \right).
\nonumber
\end{equation}
\end{enumerate}
\end{thm}

\medskip

It is worth mentioning that the particular non-diagonal form of $\vp_n(\theta)$ is, as in \cite{BroMas18}, inevitable to deduce the asymptotically non-degenerate joint distribution of the maximum-likelihood estimator (MLE), the good local maximum point $\tes$ in Theorem \ref{thm_local.LF.asymp}(2).

\begin{rem}\normalfont
\label{rem_time.scale}
Here are some comments on the model-time scale.
\begin{enumerate}
\item We are fixing the terminal sampling time $T$, so that the rate of convergence $\sqrt{n}h^{1-1/\beta}=n^{1/\beta-1/2}T^{1-1/\beta}=O(n^{1/\beta-1/2})$ for $(\lam,\mu)$.
If $\beta>1$ (resp. $\beta<1$), then a longer period would lead to a better (resp. worse) performance of estimating $(\lam,\mu)$. The Cauchy case $\beta=1$, where the two rates of convergence coincide, is exceptional.
\item We can explicitly associate a change of the terminal sampling time $T$ with those of the components of $\theta$.
Specifically, changing the model-time scale from $t$ to $tT$ in \eqref{soureg:model}, we see that the process
\begin{equation}
Y^T=(Y^T_t)_{t\in[0,1]}:=(Y_{tT})_{t\in[0,1]}
\nonumber
\end{equation}
satisfies exactly the same integral equation as in \eqref{soureg:model} except that $\theta = (\lam,\mu,\beta,\sig)$ is replaced by
\begin{equation}
\theta_T = \big( \lam_T,\mu_T,\beta_T,\sig_T):= ( T\lam,T\mu,\beta,T^{1/\beta}\sig \big)
\nonumber
\end{equation}
($\beta$ is unchanged), $X_t$ by $X^T_t:=X_{tT}$, and $J_t$ by $J^T_t := T^{-1/\beta}J_{tT}$:
\begin{equation}
Y^T_t = Y^T_0 +\int_0^t (\mu_T \cdot X^T_s-\lam_T Y^T_s)ds + \sig_T J^T_t, \qquad t\in[0,1].
\nn
\end{equation}
Note that $(J^T_t)_{t\in[0,1]}$ defines the standard $\beta$-stable {\lp}.
This indeed shows that we may set $T\equiv 1$ in the virtual (model) world without loss of generality.
This is impossible for diffusion-type models where we cannot consistently estimate the drift coefficient unless we let the terminal sampling time $T$ tend to infinity.
\end{enumerate}
\end{rem}

\begin{rem}\normalfont
The present framework allows us to do unit-period-wise, for example, day-by-day inference for both trend and scale structures, providing a sequence of period-wise estimates with theoretically valid approximate confidence sets.
This, though informally, suggests an aspect of change-point analysis in high-frequency data:
if we have high-frequency sample over $[k-1,k]$ for $k=1,\dots,[T]$, then we can construct a sequence of estimators $\{\tes(k)\}_{k=1}^{[T]}$;
then it would be possible in some way to reject the constancy of $\theta$ over $[0,[T]]$ 
if $k\mapsto \tes(k)$ ($k=1,\dots,[T]$) is not likely to stay unchanged.
\end{rem}

\begin{rem}\normalfont
It is formally straightforward to extend the model \eqref{soureg:model} to the following form:
\begin{equation}
Y_t = Y_0 +\int_0^t \left\{ a(X_s,\lam) Y_s + b(X_s,\mu) \right\}ds + \int_0^t c(X_{s-},\sig) dJ_s, \qquad t\in [0,T].
\nonumber
\end{equation}
Under mild regularity conditions on the function $(a,b,c)$ as well as on the non-random process $X$, a solution process $Y$ is explicitly given by (see \cite[Appendix]{CheKawMae03})
\begin{equation}
Y_t = e^{\psi(s,t;\lam)}Y_s + \int_s^t e^{\psi(u,t;\lam)} (b(X_u,\mu)du + c(X_{u-},\sig)dJ_u),
\qquad 0\le s<t,
\nonumber
\end{equation}
where $\psi(s,t;\lam) := \int_s^t a(X_v,\lam) dv$.
However, the corresponding likelihood asymptotics becomes much messier.
It is worth mentioning that the optimal rate matrix can be diagonal if, for example, $\p_t\p_\theta\log c(t,\sig)\not\equiv 0$ with $X_t=t$:
for details, see the previous study \cite{CleGlo20} that treated the general time-homogeneous Markovian case.
\end{rem}

\subsection{Proof of Theorem \ref{thm_local.LF.asymp}}
\label{sec_LA-proof}

In this proof, we are going to make use of the general result \cite{Swe80} about the exact-likelihood asymptotics in a more or less analogous way to that of \cite[Theorem 1]{BroMas18}:
under the uniform nature of the exact-likelihood asymptotics, we will deduce the joint convergence in distribution of the normalized score $\D_n(\tz)$ and the normalized observed information $\mci_n(\tz)$ from the uniform convergence in probability of $\mci_n(\cdot)$ in an appropriate sense.
Consequently, we will not need to derive the stable convergence in law of $\D_n(\tz)$, which is often crucial when concerned with a high-frequency sampling for a process with dependent increments.

\medskip

We have $\sup_{t\in[0,T]}|X_t|<\infty$ since $X: [0,T]\to\mbbr^q$ is assumed to be {\cadlag}.
Through the localization procedure, we may and do suppose that the driving stable {\lp} does not have jumps of size greater than some fixed threshold (see \cite[Section 6.1]{Mas19spa} for a concise account).
In that case, the L\'{e}vy measure of $J$ is compactly supported, hence in particular
\begin{equation}
\sup_{\theta\in\overline{\Theta}}\E_\theta\left[|J_1|^K\right] < \infty
\label{localization.moment}
\end{equation}
for any $K>0$.
Further, since the L\'{e}vy measure of $J$ is symmetric, the removal of large-size jumps does not change the parametric form of the drift coefficient.
We also localize the initial variable $Y_0$ so that $|Y_0|$ is essentially bounded uniformly in $\theta$.
It follows from \eqref{Y_sol} and \eqref{localization.moment} that $\sup_{\theta\in\overline{\Theta}}\sup_{0\le t\le T}\E_\theta\left[|Y_t|^K\right]<\infty$ for $t\in[0,T]$ as well.

To proceed, we introduce some further notation.
Given continuous random functions $\xi_{0}(\theta)$ and $\xi_{n}(\theta)$, $n\ge 1$, we write 
$\xi_{n}(\theta)\cip_u \xi_{0}(\theta)$ if the joint distribution of $\xi_n$ and $\xi_0$ are well-defined under $\pr_\theta$ and if $\pr_{\theta}[ |\xi_{n}(\theta)-\xi_{0}(\theta)|>\ep ] \to_{u} 0$ for every $\ep>0$ as $n\to\infty$.
Additionally, for a sequence $a_n>0$ we write $\xi_n(\theta)=o_{u,p}(a_n)$ if $a_n^{-1}\xi_{n}(\theta)\cip_u 0$, and also $\xi_n(\theta)=O_{u,p}(a_n)$ if for every $\ep>0$ there exists a constant $K>0$ for which $\sup_\theta\pr_\theta [|a_n^{-1}\xi_n(\theta)| > K]<\ep$.
Similarly, for any random functions $\chi_{nj}(\theta)$ doubly indexed by $n$ and $j\le n$, we write $\xi_{nj}(\theta)=O^\ast_p(a_n)$ if
\begin{equation}
\sup_n \max_{j\le n} \sup_{\theta} \E_\theta\left[ |a_n^{-1}\chi_{nj}(\theta)|^K\right] < \infty
\nonumber
\end{equation}
for any $K>0$. 
Finally, let
\begin{equation}
\mathfrak{N}_{n}(c;\theta) := \left\{\theta' \in \Theta:\, |\vp_{n}(\theta)^{-1}(\theta'-\theta)|\le c \right\}.
\nonumber
\end{equation}

We will complete the proof of Theorem \ref{thm_local.LF.asymp} by verifying the three statements corresponding to the conditions (12), (13), and (14) in \cite{BroMas18}, which here read
\begin{align}
& \mci_n(\theta) \cip_u \mci(\theta),
\label{C1}\\
& \sup_{\theta' \in \mathfrak{N}_{n}(c;\theta)} | \vp_{n}(\theta')^{-1}\vp_{n}(\theta) - I_{p} | \to_{u} 0,
\label{C2(i)}\\
& \sup_{\theta^1,\dots,\theta^{p}\in \mathfrak{N}_{n}(c;\theta)}
\left| \vp_{n}(\theta)^{\top}\{ \p_{\theta}^{2}\ell_{n}(\theta^{1},\dots,\theta^{p})-\p_{\theta}^{2}\ell_{n}(\theta) \} \vp_{n}(\theta) \right| 
\v2{\cip_{u}} 0,
\label{C2(ii)}
\end{align}
respectively, where \eqref{C2(i)} and \eqref{C2(ii)} should hold for all $c>0$ and where $\p_{\theta}^{2}\ell_{n}(\theta^1,\dots,\theta^{p})$, $\theta^k\in\Theta$, denotes the $p\times p$ Hessian matrix of $\ell_n(\theta)$, whose $(k,l)$th element is given by $\p_{\theta_k}\p_{\theta_l}\ell_{n}(\theta^k)$, where $\theta=:(\theta_l)_{l=1}^{p}$.
Having obtained \eqref{C1}, \eqref{C2(i)} and \eqref{C2(ii)}, \cite[Theorem 1 and 2]{Swe80} immediately concludes Theorem \ref{thm_local.LF.asymp}.
We can verify \eqref{C2(i)} exactly as in \cite{BroMas18}, so will look at \eqref{C1} and \eqref{C2(ii)}.

\medskip

\noindent
\textit{Proof of \eqref{C1}.}
Recall the expression \eqref{log-lf}.
To look at the entries of $\p_\theta^2\ell_n(\theta)$, we introduce several shorthands for notational convenience; they may look somewhat daring, but 
would not bring confusion.
Let us omit the subscript $\beta$ and the argument $\ep_j$ of the aforementioned notation, such as $\phi:=\phi_\beta(\ep_j)$, $g:=g_\beta(\ep_j)$ and so on. For brevity, we also write 
\begin{equation}
l'=\log(1/h), \quad 
\v2{c=\eta(\lam\beta h)^{-1/\beta}}, \quad \ep=\ep_{j}(\theta),
\nonumber
\end{equation}
so that \eqref{log-lf} becomes
\begin{equation}
\ell_n(\theta) = \sumj \left( -\log\sig +\frac1\beta l' + \log c + \log\phi \right).
\nonumber
\end{equation}
Further, the partial differentiation with respect to a variable will be denoted by the braced subscript such as $\ep_{(a)}:=\p_a\ep_j(\theta)$ and $\ep_{(a,b)}:=\p_a\p_b\ep_j(\theta)$.
Then, direct computations give the first-order partial derivatives:
\begin{align}
\v2{\p_\lam}\ell_n(\theta) &= \sumj \left( (\log c)_{(\lam)} + \ep_{(\lam)}\, g \right), \nn\\
\v2{\p_\mu}\ell_n(\theta) &= \sumj \ep_{(\mu)}\, g, \nn\\
\p_\beta\ell_n(\theta) &= \sumj \left( -\beta^{-2}l' + (\log c)_{(\beta)} + \ep_{(\beta)}\, g + f \right), \nn\\
\p_\sig\ell_n(\theta) &= \sumj \left( -\sig^{-1} + \ep_{(\sig)} \, g \right),
\nonumber
\end{align}
followed by the second-order ones:
\begin{align}
\p_\lam^2\ell_n(\theta) &= \sumj \left\{ (\ep_{(\lam)})^2 (\p g) + \ep_{(\lam,\lam)}\, g + (\log c)_{(\lam,\lam)} \right\},
\nn
\\
\p_\mu^2\ell_n(\theta) &= \sumj (\ep_{(\mu)})^2 (\p g),
\nn
\\
\p_\beta^2\ell_n(\theta) &= \sumj \left\{ 
2\v2{\beta^{-3}} l' + (\log c)_{(\beta,\beta)} + \ep_{(\beta,\beta)}\, g + \ep_{(\beta)}\, g_{(\beta)}
+ (\ep_{(\beta)})^2 (\p g) + f_{(\beta)} + \ep_{(\beta)}\, (\p f) \right\},
\nn
\\
\p_\sig^2\ell_n(\theta) &= \sumj \left\{ \sig^{-2} + (\ep_{(\sig)})^2 (\p g) + \ep_{(\sig,\sig)} \, g \right\},
\nn
\\
\p_\lam\p_\mu\ell_n(\theta) &= \sumj \left\{ \ep_{(\lam)}\,\ep_{(\mu)} (\p g) + \ep_{(\mu,\lam)}\, g \right\},
\nn
\\
\p_\lam\p_\beta\ell_n(\theta) &= \sumj \left\{
(\log c)_{(\lam,\beta)} + \ep_{(\lam)}\, g_{(\beta)} + \ep_{(\beta)}\, \ep_{(\lam)}\, (\p g) + \ep_{(\lam,\beta)}\, g
\right\},
\nn
\\
\p_\lam\p_\sig\ell_n(\theta) &= \sumj \left\{ \ep_{(\sig)}\,\ep_{(\lam)} (\p g) + \ep_{(\lam,\sig)}\, g \right\},
\nn
\\
\p_\mu\p_\beta\ell_n(\theta) &= \sumj \left\{ \ep_{(\mu)}\, g_{(\beta)} + \ep_{(\beta)}\,\ep_{(\mu)} (\p g) + \ep_{(\beta,\mu)}\, g \right\},
\nn
\\
\p_\mu\p_\sig\ell_n(\theta) &= \sumj \left\{ \ep_{(\mu)}\,\ep_{(\sig)} (\p g) + \ep_{(\mu,\sig)} \, g \right\},
\nn
\\
\p_\beta\p_\sig\ell_n(\theta) &= \sumj \left\{ \ep_{(\sig)}\, g_{(\beta)} + \ep_{(\beta)}\,\ep_{(\sig)} (\p g) 
+ \ep_{(\beta,\sig)}\, g \right\}.
\nn
\end{align}
It is straightforward to see which term is the leading one in each expression above.
We do not list all the details here, but for later reference mention a few of the points:

\begin{itemize}

\item $\p^k\log\eta(y) = O_u(1)$ for $|y|\to 0$ whatever $k\in\mbbzp$ is;

\item $(\log c)_{(\lam,\dots,\lam)}=O_u(h^k)$ ($k$-times, $k\in\mbbzp$), $(\log c)_{(\lam,\beta)}=O_u(h^2)$, $(\log c)_{(\beta)}=O_u(h^2)$, $(\log c)_{(\beta,\beta)}=O_u(h^4)$, and so forth;

\item $\max_{j\le n}|\p_\lam^k\zeta_j(\lam)|=O(h^k)$ for $k\in\mbbzp$;

\item \hmcheck
Recalling the definition \eqref{epj_def} and because of the consequence \eqref{localization.moment} of the localization, concerning the partial derivatives of $\ep_j(\theta)$ we obtain the asymptotic representations:
$\ep_{(\mu,\sig)} = (1+o_{u}(1)) \sig^{-2}h^{1-1/\beta}$,
$\ep_{(\mu,\lam)} = (1+o_{u}(1)) \sig^{-1}h^{2-1/\beta}/2$,
$\ep_{(\sig,\lam)} = (1+o_{u}(1))
\{-\sig^{-2}h^{1-1/\beta}Y_{t_{j-1}} + O^\ast_p(h \vee h^{2-1/\beta})\}$,
$\ep_{(\lam,\lam)} = O^\ast_p(h^{2-1/\beta})$,
$\ep_{(\beta,\beta)} = O_p^\ast(h^2 (l')^2) + \ep\, O^\ast_p((l')^2)$, 
$\ep_{(\lam,\beta)} = O^\ast_p(l' h^{1-1/\beta})$, and so on;
the terms ``$o_u(1)$'' therein are all valid uniformly in $j\le n$.
\end{itemize}

Now we write
\begin{equation}
\mci_n(\theta) = 
\begin{pmatrix}
\mci_{11,n}(\theta) & \mci_{12,n}(\theta) \\
\mci_{12,n}(\theta)^\top & \mci_{22,n}(\theta)
\end{pmatrix}
\nonumber
\end{equation}
with $\mci_{11,n}(\theta) \in \mbbr^{1+q}\otimes\mbbr^{1+q}$, $\mci_{22,n}(\theta) \in \mbbr^{2}\otimes\mbbr^{2}$ and $\mci_{12,n}(\theta) \in \mbbr^{1+q}\otimes\mbbr^{2}$.
We can deduce $\mci_{22,n}(\theta) \cip_u \mci_{\beta,\sig}(\theta)$ in exactly the same way as in the proof of Eq.(12) in \cite{BroMas18}.
Below, we will show $\mci_{11,n}(\theta) \cip_u \mci_{\lam,\mu}(\theta)$ and $\mci_{12,n}(\theta) \cip_u 0$.

\medskip

The Burkholder inequality ensures that
\begin{equation}
\frac{1}{\sqrt{n}}\sumj \pi(X_{t_{j-1}},Y_{t_{j-1}};\theta)U(\ep_j(\theta)) = O_{u,p}(1)
\label{add-1}
\end{equation}
for any continuous $\pi(x,y;\theta)$ and for any $U(\ep_j(\theta))$ such that $\E_\theta[U(\ep_j(\theta))]=0$ ($\theta\in\Theta$) and that the left-hand side of \eqref{add-1} is continuous over $\theta\in\overline{\Theta}$. 
Also, note that the right continuity of $t\mapsto X_{t}$ implies that ($X_t^{\otimes 1}:=X_t$)
\begin{equation}
\lim_{n\to\infty}\max_{l=1,2}\max_{j\le n}\left|\frac1h \int_j (X_s^{\otimes l} - X_{t_{j-1}}^{\otimes l}) ds\right| = 0.
\nn
\end{equation}
These basic facts will be repeatedly used below without mentioning them.

For convenience, we will write
\begin{equation}
r_n=r_n(\beta)=\sqrt{n}h^{1-1/\beta}
\label{def:r_n}
\end{equation}
and denote by $\bm{1}_{u,p}$ any random array $\xi_{nj}(\theta)$ such that $\max_{j\le n}|\xi_{nj}(\theta) - 1| \cip_u 0$.
Direct computations give the following expressions for the components of $\mci_{11,n}(\theta)=- r_{n}^{-2}\p_{(\lam,\mu)}^2\ell_{n}(\theta)$:
\begin{align}
-\frac{1}{r_n^2}\p_\mu^2\ell_n(\theta)
&= -\frac1n\sumj \sig^{-2}(\p g)\,\zeta_j(\lam)^{\otimes 2} + o_{u,p}(1) \nn\\
&= \frac1n\sumj \sig^{-2} g^2 \,\zeta_j(\lam)^{\otimes 2} + o_{u,p}(1),
\nn\\
-\frac{1}{r_n^2}\p_\lam^2\ell_n(\theta) 
&= -\frac1n\sumj \sig^{-2}(\p g)\,Y_{t_{j-1}}^2\,\bm{1}_{u,p} + o_{u,p}(1)+O_{u,p}(h^{1/\beta}) \label{ell.p2.trend_1} \nn\\
&= \frac1n\sumj \sig^{-2} g^2\, Y_{t_{j-1}}^2 + o_{u,p}(1),
\nn\\
-\frac{1}{r_n^2}\p_\lam\p_\mu\ell_n(\theta) 
&= -\frac1n\sumj \left\{
(\p g) \left( \bm{1}_{u,p}\sig^{-1}Y_{t_{j-1}}\,\zeta_j(\lam) + O^\ast_{p}(h^{1/\beta}) \right) (-\sig^{-1}\bm{1}_{u,p})
+ O^\ast_{p}(h^{1/\beta}) \right\} \nn\\
&=\bm{1}_{u,p}\left( - \frac1n\sumj \sig^{-2}(\p g)\,Y_{t_{j-1}}\,\zeta_j(\lam) + o_{u,p}(1) \right) + o_{u,p}(1) \nn\\
&= \frac1n\sumj \sig^{-2} g^2\, Y_{t_{j-1}}\,\zeta_j(\lam) + o_{u,p}(1).
\nn
\end{align}
We can deduce that $\mci_{11,n}(\theta) \cip_u \mci_{\lam,\mu}(\theta)$ as follows.
\begin{itemize}
\item First, noting that $\ep_j=\ep_j(\theta) \stackrel{\pr_{\theta}}{\sim} \text{i.i.d.}~\mcl(J_1)$, we make the compensation $g^2= \E_\theta[g^2]+(g^2 - \E_\theta[g^2])$ in the summands in rightmost sides of the last three displays and then pick up the leading part involving $\E_\theta[g^2]$;
the other one becomes negligible by the Burkholder inequality.
\item \hmcheck
Then, 
\v2{the a.s. Riemann integrability of $t\mapsto (X_{t}(\omega),Y_{t}(\omega))$ allows us to conclude that, for $k,l\in\{0,1,2\}$ and under $\pr_\theta$ for each $\theta$,}
\begin{align}
D_{n}(k,l)&:= \left| \frac1n \sumj Y_{t_{j-1}}^k\,X_{t_{j-1}}^{\otimes l} - \frac1T \int_0^T Y_{t}^k \, X_{t}^{\otimes l} dt \right|
\nn\\
&\lesssim 
\frac1n \sumj \frac1h \int_j  \Bigg(|Y_t- Y_{t_{j-1}}|(1+|Y_t|+|Y_{t_{j-1}}|)^C \nn\\
&{}\qquad + |Y_t|^k \left|\left(\frac1h\int_j X_t dt + O(h)\right)^{\otimes l} - X_t ^{\otimes l}
\right|\Bigg)dt
\nn\\
&\lesssim 
\frac1n \sumj \frac1h \int_j  \Bigg(|Y_t- Y_{t_{j-1}}|(1+|Y_t|+|Y_{t_{j-1}}|)^C + |Y_t|^k o(1)\Bigg)dt
\cip 0,
\nn
\end{align}
where the order symbols in the estimates are valid uniformly in $j\le n$.
By \eqref{Y_sol}, under the localization we have $\max_{j\le n}\sup_\theta \E_\theta[|Y_t|^M]=O(1)$ and $\max_{j\le n}\sup_\theta \E_\theta[|Y_t- Y_{t_{j-1}}|^M]=o_u(1)$ for any $M>0$, from which it follows that $D_{n}(k,l) = o_{u,p}(1)$.
\end{itemize}
Specifically, for the case of $-r_n^{-2}\p_\mu^2\ell_n(\theta)$, we have
\begin{align}
-\frac{1}{r_n^2}\p_\mu^2\ell_n(\theta)
&= \frac1n\sumj \sig^{-2} \E_\theta[g^2] \,\zeta_j(\lam)^{\otimes 2} + o_{u,p}(1)
\nn\\
&= \sig^{-2} \E_\theta[g_\beta(\ep_1(\theta))^2] \, \frac{1}{T}\int_0^T X_{t}^{\otimes 2}dt + o_{u,p}(1)
\cip_u \mci_{\lam,\mu;22}(\theta)
\nonumber
\end{align}
with $\mci_{\lam,\mu;22}(\theta)$ denoting the lower left $q\times q$ component of $\mci_{\lam,\mu}(\theta)$.
The others can be handled analogously.

Next we turn to looking at $\mci_{12,n}(\theta)=\{\mci_{12,n}^{kl}(\theta)\}_{k,l}$:
\begin{align}
\mci_{12,n}^{11}(\theta) &= \vp_{11,n}(\theta)\p_\lam\p_\beta\ell_n(\theta) + \vp_{21,n}(\theta)\p_\mu\p_\beta\ell_n(\theta), \nn\\
\mci_{12,n}^{12}(\theta) &= \vp_{11,n}(\theta)\p_\lam\p_\sig\ell_n(\theta) + \vp_{21,n}(\theta)\p_\mu\p_\sig\ell_n(\theta), \nn\\
\mci_{12,n}^{21}(\theta) &= \vp_{12,n}(\theta)\p_\lam\p_\beta\ell_n(\theta) + \vp_{22,n}(\theta)\p_\mu\p_\beta\ell_n(\theta), \nn\\
\mci_{12,n}^{22}(\theta) &= \vp_{12,n}(\theta)\p_\lam\p_\sig\ell_n(\theta) + \vp_{22,n}(\theta)\p_\mu\p_\sig\ell_n(\theta).
\nonumber
\end{align}
We can deduce that $\mci_{12,n}(\theta) \cip_u 0$ just by inspecting the four components separately in a similar way that we managed $\mci_{11,n}(\theta)$.
Let us only mention the lower-left $q\times 1$ component: recalling the properties \eqref{hm:dens.ab} and $|\vp_{22,n}|\lesssim_u l'$, we see that
\begin{align}
\mci_{12,n}^{21}(\theta)
&= -\frac{(h^{1-1/\beta})^{-1}}{n}\sumj \bigg( \vp_{12,n}\,(\log c)_{(\lam,\beta)} + \vp_{12,n}\,\ep_{(\lam)}\,g_{(\beta)}
+ \vp_{12,n}\,\ep_{(\lam)}\,\ep_{(\beta)}\,(\p g) \nn\\
&{}\qquad + \vp_{12,n}\, \ep_{(\lam,\beta)}\, g
+ \vp_{22,n}\,\ep_{(\mu)}\, g_{(\beta)} + \vp_{22,n}\,\ep_{(\beta)}\,\ep_{(\mu)}\,(\p g) + \vp_{22,n}\,\ep_{(\beta,\mu)} \,g\bigg) \nn\\
&= O(h^{1+1/\beta}) + O_{u,p}(n^{-1/2}\vee h^{1/\beta}) + O_{u,p}\left( (n^{-1/2}\vee h^{1/\beta}) \,l'\right) \nn\\
&{}\qquad +O_{u,p}(n^{-1/2}) + O_{u,p}\left(n^{-1/2}(l')^2\right) + O_{u,p}\left(n^{-1/2}(l')^2\right) \cip_u 0.
\nonumber
\end{align}
Thus, the claim \eqref{C1} follows.

\medskip

\noindent
\textit{Proof of \eqref{C2(ii)}.}
\hmcheck
Note that
\begin{equation}
\sup_{\theta'\in\mathfrak{N}_{n}(c;\theta)}
|\ep_{j}(\theta')| \lesssim_{u} |\ep_{j}(\theta)| + \overline{s}_{nj}(\theta;c),
\nonumber
\end{equation}
where $|\overline{s}_{nj}(\theta;c)| \lesssim o_{u}(1)(1+|Y_{t_{j-1}}|)$.
Also, for each $k,l,m\in\mbbzp$, we have $\pr_\theta$-a.s. the (rough) estimate:
\begin{equation}
\frac{1}{n}\big|\p_{\beta}^{k}\p_{\sig}^{l}\p_{(\lam,\mu)}^{m}\ell_{n}(\theta)\big| 
\lesssim_{u} (l')^{k}\, h^{(1-1/\beta)m}\,
\frac{1}{n}\sumj (1+|Y_{t_{j-1}}|)^m \left\{1 + \log\left(1+|\ep_{j}(\theta)|^{2}\right)\right\}^{k}.
\nn
\end{equation}
Then, as in the proof of Eq.(14) in \cite{BroMas18}, for each $c>0$ we can find a constant $R=R(c)>0$ such that (still rough, but sufficient)
\begin{align}
& \sup_{\theta^1,\dots,\theta^{p}\in \mathfrak{N}_{n}(c;\theta)}
\left| \vp_{n}(\theta)^{\top}\{ \p_{\theta}^{2}\ell_{n}(\theta^{1},\dots,\theta^{p})-\p_{\theta}^{2}\ell_{n}(\theta) \} \vp_{n}(\theta) \right|
\nn\\
& \lesssim_{u} 
\sup_{\theta',\theta^1,\dots,\theta^{p}\in \mathfrak{N}_{n}(c;\theta)}
\left| \vp_{n}(\theta)^{\top} \left\{ \p_{\theta}^{3}\ell_{n}(\theta^{1},\dots,\theta^{p}) [\theta' - \theta] \right\} \vp_{n}(\theta) \right|
\nn\\
&\lesssim_{u} \frac{(l')^{C}}{\sqrt{n}} 
\sup_{\beta',\beta'' \in \overline{B}(\beta;R/l^{\prime})} h^{(1/\beta'-1/\beta'')3} 
\sup_{\theta' \in \mathfrak{N}_{n}(c;\theta)} \frac{1}{n}\sumj (1+|Y_{t_{j-1}}|)^m \left\{1 + \log\left(1+|\ep_{j}(\theta')|\right)\right\}^{3} \nn\\
&\lesssim_{u} \frac{(l')^{C}}{\sqrt{n}} \frac{1}{n}\sumj (1+|Y_{t_{j-1}}|)^m \left\{1 + \log\left(1+|\ep_{j}(\theta)|\right)\right\}^{3}
\lesssim O_{u,p}\left(\frac{(l')^{C}}{\sqrt{n}}\right) \cip_u 0,
\nn
\end{align}
where $\overline{B}(\beta;R/l_n^{\prime})$ denotes the closed ball with center $\beta$ and radius $R/l'$.
This shows \eqref{C2(ii)}. The proof of Theorem \ref{thm_local.LF.asymp} is complete.

\section{Asymptotically efficient estimator}
\label{sec:efficient.estimator}

From now on, we fix a true value $\tz\in\Theta$, and the stochastic symbols and convergences will be taken under $\pr:=\pr_{\tz}$; accordingly, we write $\E:=\E_{\tz}$.
Having Theorem \ref{thm_local.LF.asymp} in hand, we can proceed with the construction of an asymptotically efficient estimator.
It is known that any asymptotically centering estimator $\tes^\ast$:
\begin{equation}
\vp_{n}(\tz)^{-1}(\tes^\ast -\tz) = \mci_{n}(\tz)^{-1}\D_n(\tz) + o_{p}(1)
\label{def:ACe}
\end{equation}
are regular; 
by Theorem \ref{thm_local.LF.asymp}, the right-hand side converges in distribution to $MN_{p,\tz}\left(0,\, \mci(\tz)^{-1} \right)$.
This together with the convolution theorem in turn gives the asymptotic minimax theorem:
for any measurable (loss) function $\mfL:\,\mbbr^p \to \mbbrp$ such that $\mfL(u)=\tau(|u|)$ for some non-decreasing $\tau:\,\mbbrp\to\mbbrp$ with $\tau(0)=0$, we have
\begin{equation}
\liminf_{n\to\infty} \E\left[ \mfL\left( \vp_{n}(\tz)^{-1}(\tes^\ast - \tz) \right)\right]
\ge \E\big[\mfL\big(\mci(\tz)^{-1/2}Z\big)\big].
\label{lower.bound}
\end{equation}
\v2{
Recalling that 
$\mcl\left( \D_{n}(\theta), \, \mci_n(\theta) |\pr_\theta\right) \wc \mcl\left( \mci(\theta)^{1/2}Z,\, \mci(\theta) \right)$,
where $Z\sim N_{p}(0,I)$ (Theorem \ref{thm_local.LF.asymp}) and in view of the lower bound in \eqref{lower.bound},
we may call that any estimator $\tes^\ast$ satisfying \eqref{def:ACe} asymptotically efficient.}
Again by Theorem \ref{thm_local.LF.asymp}, the good local maximum point $\tes$ of $\ell_n(\theta)$ is asymptotically efficient.
We refer to \cite[Theorems 2 and 3, and Proposition 2]{Jeg82} and also \cite[Theorem 8]{Jeg95} for more information and details of the above arguments.

Theorem \ref{thm_local.LF.asymp} is based on the classical Cram\'{e}r-type argument. The well-known shortcoming is its local character: the result just tells us the existence of an asymptotically nicely behaving root of the likelihood equation, but does not give information about which local maxima is the one when there are multiple local maxima, equivalently multiple roots for the likelihood equations \cite[Section 7.3]{Leh99}.
Indeed, the log-likelihood function $\ell_n$ of \eqref{log-lf} is highly nonlinear and non-concave.
In this section, we try to get rid of the locality by a Newton-Raphson-type improvement, which in our case will not only remedy the aforementioned inconvenience of the multiple-root problem but also enable us to bypass the numerical optimization involving the stable density $\phi_\beta$.
In \cite[Section 3]{BroMas18}, for the $\beta$-stable {\lp} (the special case of \eqref{soureg:model} with $\lam=0$ and $X\equiv 1$), we provided an initial estimator based on the sample median and the method of moments associated with logarithm and/or lower-order fractional moments.
However, it was essential in \cite{BroMas18} that the model is a {\lp} for which we could apply the median-adjusted central limit theorem for an i.i.d. sequence of random variables.
In the present case, we need a 
different sort of argument.

In Theorem \ref{thm_local.LF.asymp}, the process $X=(X_t)_{t\in[0,T]}$ was assumed to be observed continuously in $[0,T]$.
In this section, we will instead deal with a discrete-time sample $(X_{t_j})_{j=0}^{n}$ under the additional condition:
\begin{equation}
\exists \kappa \in \v2{(1/2,1]},\quad 
\max_{j\le n}\left| \frac1h \int_j (X_t - X_{t_{j-1}}) dt \right|
\lesssim h^\kappa.
\label{X-kappa}
\end{equation}
We will explicitly construct an estimator $\tes^\ast$ which is asymptotically equivalent to the MLE $\tes$, by verifying the asymptotically centering property \eqref{def:ACe}; for this much-thinned sample, we may and do keep calling such a $\tes^\ast$ asymptotically efficient.

\subsection{Newton-Raphson procedure}

To proceed with a discrete-time sample $\{(X_{t_j},Y_{t_j})\}_{j=0}^{n}$, we introduce the approximate-likelihood function $\mbbh_n(\theta)$ by replacing $\zeta_{j}(\lam)$ by $X_{t_{j-1}}$ in the definition \eqref{log-lf} of the genuine log-likelihood function $\ell_n(\theta)$ (recall the notation $l':=\log(1/h)$):
\begin{align}
\mbbh_n(\theta) 
&=\sumj \left( -\log\sig +\frac1\beta l' 
\v2{- \frac1\beta \log\eta(\lam\beta h)} + \log\phi_{\beta} \left( \ep_{j}'(\theta) \right)\right),
\label{log-qlf}
\end{align}
where
\begin{equation}
\ep'_{j}(\theta) := \frac{Y_{t_{j}} - e^{-\lam h}Y_{t_{j-1}} - \mu \cdot X_{t_{j-1}} h}{\sig h^{1/\beta}\eta(\lam\beta h)^{1/\beta}}.
\label{epj'_def}
\end{equation}
Of course, this approximation is not for free: to manage the resulting discretization error specified later on, we additionally impose that
\begin{equation}
\beta_0 > \frac{2}{1+2\kappa}.
\label{beta&kappa}
\end{equation}
Then we have at least $\beta_0 > 2/3$, so that small-values of $\beta_0$ are excluded;
this is the price we have to pay for dealing with a discrete-time sample from $X$ in an efficient way.
Accordingly, in the sequel, we will reset the parameter space of $\beta$ to be a domain $\Theta_\beta$ such that $\overline{\Theta_\beta} \subset(2/3,2)$.

Toward construction of an asymptotically efficient estimator $\tes^\ast$ satisfying \eqref{def:ACe}, we will prove a basic result about a Newton-Raphson type procedure. As in \eqref{def:r_n}, we write $r_n=r_n(\beta_0)=\sqrt{n}h^{1-1/\beta_0}$. Write $n^{-1/2}\tilde{\vp}_n(\theta)$ for the lower-right $2\times 2$-part of $\vp_n(\theta)$,
so that the definition \eqref{hm:vp_def} with $\theta=\tz$ becomes $\vp_n(\tz) = \diag(r_{n}^{-1}I_{q+1},\, n^{-1/2}\tilde{\vp}_n(\tz))$.
We then introduce the diagonal matrix
\begin{equation}
\vp_{0,n}=\vp_{0,n}(\beta_0):= \diag\left(
r_{n}^{-1}I_{q+1},\, n^{-r/2}
\begin{pmatrix}
1 & 0 \\ 0 & l'
\end{pmatrix}
\right)
\label{vp0.rate}
\end{equation}
for a constant
\begin{equation}
\v2{0<r\le 1.}
\label{r_set}
\end{equation}
The difference between $\vp_n$ and $\vp_{0,n}$ is only in the lower-right component for $(\beta,\sig)$, and note that the matrix $\vp_n^{-1}\vp_{0,n}$ may diverge in norm.
Then, suppose that we are given an initial estimator $\tesnn_{0,n}=(\lesnn_{0,n},\mesnn_{0,n},\besnn_{0,n},\sesnn_{0,n})$ such that $\vp_{0,n}^{-1}(\tesnn_{0,n} - \tz) = O_p(1)$, namely
\begin{equation}
\left(r_n (\lesnn_{0,n}-\lam_0),\, r_n(\mesnn_{0,n}-\mu_0),\,
n^{r/2}(\besnn_{0,n}-\beta_0),\,\frac{n^{r/2}}{l'}(\sesnn_{0,n}-\sig_0) \right) = O_{p}(1).
\nn
\end{equation}

\v2{Let us write $a=(\lam,\mu)$ and $b=(\beta,\sig)$.}
Based on the approximate-likelihood function \eqref{log-qlf} and $\tesnn_{0,n}$, we recursively define the $k$-step estimator $\tesnn_{k,n}$ ($k\ge 1$) by
\v2{
\begin{align}
\tesnn_{k,n} = \tesnn_{k-1,n} + \left\{
\diag\left(-\p_a^2\mbbh_n(\tesnn_{k-1,n}),\, -\p_b^2\mbbh_n(\tesnn_{k-1,n})
\right)\right\}^{-1} \p_\theta\mbbh_n(\tesnn_{k-1,n})
\label{tes.k_def}
\end{align}
on the event 
$F_{k-1,n} := \{|\det(\p_a^2\mbbh_n(\tesnn_{k-1,n}))| \wedge |\det(\p_b^2\mbbh_n(\tesnn_{k-1,n}))| > 0\}$
} 
and assign an arbitrary value to $\tesnn_{k,n}$ on the complement set $F_{k-1,n}^c$; 
below, it will be seen (as in the proof of Theorem \ref{thm_local.LF.asymp}) that $\pr[F_{k-1,n}]\to 1$, hence the arbitrary property does not matter asymptotically and we may and do suppose that $\pr[F_{k-1,n}]=1$ for $k\ge 1$.
%
\v2{
In our subsequent arguments, the inverse-matrix part in \eqref{tes.k_def} must be block-diagonal:
see Remark \ref{rem_block.diag.form} below.
}

In what follows, $\tes$ denotes the good local maxima of the likelihood function $\ell_n(\theta)$, when $(X_t)_{t\le T}$ is observable; by Theorem \ref{thm_local.LF.asymp}, we have $\pr[\p_\theta\ell_n(\tes)=0] \to 1$ and $\vp_{n}(\tz)^{-1}(\tes -\tz) = \mci_{n}(\tz)^{-1}\D_n(\tz) + o_{p}(1)\wc MN_{p,\tz}\left(0,\, \mci(\tz)^{-1} \right)$.
\v2{
Define the number
\begin{equation}
K := \min\{ k\in\mbbn :\, 2^{k-1}r > 1/2\} = \min\{ k\in\mbbn :\, k>\log_2(1/r)\}.
\label{def_K}
\end{equation}
We are going to deduce the asymptotic equivalence of $\tes$ and $\tesnn_{K,n}$:
\begin{equation}
\vp_{n}(\tz)^{-1}(\tesnn_{K,n} -\tes) = o_{p}(1),
\label{NR.proof-5}
\end{equation}
starting from the initial estimator $\tesnn_{0,n}$;
\eqref{NR.proof-5} concludes \eqref{def:ACe} (hence \eqref{lower.bound} as well) with \v2{$\tes^\ast = \tesnn_{K,n}$}.
}


\v2{
We are assuming that $\vp_{0,n}^{-1}(\tesnn_{0,n}-\tz)=O_p(1)$ (hence also $\vp_{0,n}^{-1}(\tesnn_{0,n}-\tes)=O_p(1)$).
Then, to establish \eqref{NR.proof-5}, we first look at the amount of improvement through \eqref{tes.k_def} with $k=1$.
Write $\vp_n = \vp_{n}(\tz)$ and $\tilde{\vp}_n = \tilde{\vp}_{n}(\tz)$, and introduce
\begin{align}
\hat{\mci}_{0,n} &:= -\vp_{n}^{\top}
\diag\left( \p_a^2\mbbh_n(\tesnn_{0,n}),\, \p_b^2\mbbh_n(\tesnn_{0,n}) \right)\vp_{n}
\nn\\
&= \diag\left( - 
r_n^{-2}
\p_a^2\mbbh_n(\tesnn_{0,n}),\, 
-\frac{1}{n} \tilde{\vp}_n^{\top} \p_b^2\mbbh_n(\tesnn_{0,n}) \tilde{\vp}_n \right)
=: \big( \hat{\mci}_{0,a,n},\, \hat{\mci}_{0,b,n} \big).
\nonumber
\end{align}
We apply Taylor's expansion around $\tes$ to \eqref{tes.k_def} with $k=1$:
for some random point $\tesnn'_{0,n}$ on the segment joining $\tesnn_{0,n}$ and $\tes$,
\begin{align}
\tesnn_{1,n} -\tes 
&= \tesnn_{0,n} -\tes + \vp_n \hat{\mci}_{0,n}^{-1} \vp_n^\top \, \p_\theta\mbbh_n(\tesnn_{0,n})
\nn\\
&= \vp_n \hat{\mci}_{0,n}^{-1}
\Big\{
\vp_n^\top \p_\theta\mbbh_n(\tes)
\nn\\
&{}\qquad + \vp_n^\top\left( 
\diag\left( -\p_a^2\mbbh_n(\tesnn_{0,n}),\, -\p_b^2\mbbh_n(\tesnn_{0,n}) \right)
- \big( - \p_{\theta}^{2}\mbbh_n(\tesnn'_{0,n}) \big) 
\right) [\tesnn_{0,n}-\tes]
\Big\}
\nn\\
&=: \vp_n \hat{\mci}_{0,n}^{-1}\big( R'_{0,n} + R''_{0,n} \big).
\label{NR.proof-1}
\end{align}
In what follows, we will derive the rate of convergence of $\tesnn_{1,n} -\tes$ in several steps.
Here again, we may and do work under the localization (See Section \ref{sec_LA-proof}).
}

\medskip

\v2{
\paragraph{\textit{Step 1.}}
First, we show that $\hat{\mci}_{0,n}^{-1}=O_p(1)$. We have
\begin{align}
\hat{\mci}_{0,a,n} &= - r_n^{-2} \p_a^2\mbbh_n(\tz) - r_n^{-2} \left(\p_a^2\mbbh_n(\tesnn_{0,n}) - \p_a^2\mbbh_n(\tz)\right),
\label{v2-6}\\
\hat{\mci}_{0,b,n} &= - \frac{1}{n} \tilde{\vp}_n^{\top}\p_b^2\mbbh_n(\tz) \tilde{\vp}_n 
- \frac{1}{n} \tilde{\vp}_n^{\top} \left(\p_b^2\mbbh_n(\tesnn_{0,n}) - \p_b^2\mbbh_n(\tz) \right) \tilde{\vp}_n.
\label{v2-7}
\end{align}
The first terms on the right-hand sides above tend to $\mci_{\lam,\mu}(\tz)$ and $\mci_{\beta,\sig}(\tz)$ in probability, respectively.
\hmcheck The second terms equal $o_p(1)$, by similar considerations to the verification of \eqref{C2(ii)} in the proof of Theorem \ref{thm_local.LF.asymp}.
Hence $\hat{\mci}_{0,n} \cip \mci(\tz)$ and in particular $\hat{\mci}_{0,n}^{-1}=O_p(1)$ since $\mci(\tz)$ is a.s. positive definite.
}

\medskip

\v2{
\paragraph{\textit{Step 2.}}
Next, we show that $R'_{0,n} =\vp_n^\top \p_\theta\mbbh_n(\tes))$ is $o_p(1)$. Observe that
\begin{align}
\vp_n^\top \p_\theta\mbbh_n(\tes) &= \vp_n^\top \p_\theta\ell_n(\tes) + \vp_n^\top \left( \p_\theta\mbbh_n(\tes) - \p_\theta\ell_n(\tes)\right).
\nn
\end{align}
For the first term, we have $\vp_n^\top \p_\theta\ell_n(\tes) = o_p(1)$ since $\pr[|s_n \p_\theta\ell_n(\tes)|>\ep] \le \pr[|\p_\theta\ell_n(\tes)|\ne 0]\to 0$ for every $\ep>0$ and $s_n\uparrow\infty$. 
To manage the second term, we need to estimate the gap between $\mbbh_n(\theta)$ and $\ell_n(\theta)$ by taking the different convergence rates of their components into account.
By the definitions \eqref{log-lf} and \eqref{log-qlf},
\begin{align}
\mbbh_n(\theta) - \ell_n(\theta) &= \sumj \left( \log\phi_\beta(\ep'_j(\theta)) - \log\phi_\beta(\ep_j(\theta)) \right)
\nn\\
&= \sumj \left(\int_0^1 g_\beta\left( \ep_j(\theta) + s (\ep'_j(\theta)-\ep_j(\theta))\right)ds\right) (\ep'_j(\theta)-\ep_j(\theta).)
\label{v2-3}
\end{align}
From the expressions \eqref{epj_def} and \eqref{epj'_def} and since $\kappa\le 1$, a series of straightforward computations show that the partial derivatives of
\begin{align}
d_{\ep,j}(\theta)
&:=\ep'_j(\theta)-\ep_j(\theta)
\nn\\
&=
\frac{1}{\sig h^{1/\beta}\eta(\lam\beta h)^{1/\beta}}
\left(\mu \cdot \int_j(e^{-\lam(t_j -s)} -1)X_s ds + \mu\cdot \int_j (X_s - X_{t_{j-1}})ds\right)
\nonumber
\end{align}
satisfies the following bounds:
$|\p_\mu d_{\ep,j}(\theta)| \lesssim h^{1+\kappa-1/\beta}$, 
$|\p_\lam d_{\ep,j}(\theta)| \lesssim h^{2-1/\beta}$, 
$|\p_\beta d_{\ep,j}(\theta)| \lesssim h^{1+\kappa-1/\beta} l'$, 
and $|\p_\sig d_{\ep,j}(\theta)| \lesssim h^{1+\kappa-1/\beta}$.
Obviously, $h^{1/\tilde{\beta}_n - 1/\beta_0} = 1 + o_p(1)$ for any $\tilde{\beta}_n$ such that $n^{v}(\tilde{\beta}_n-\beta_0)=O_p(1)$ for some $v>0$;
below, we will repeatedly make use of this fact without mention.
Further, under \eqref{beta&kappa}, it holds that
\begin{equation}
\exists \del_1>0,\quad \sqrt{n} \,h^{1+\kappa -1/\beta_0} = O(n^{-1/2-\kappa +1/\beta_0}) =O(n^{-\del_1}).
\label{v2-1}
\end{equation}
By piecing together these observations, the basic property \eqref{hm:dens.ab}, and the expression \eqref{v2-3}, under \eqref{X-kappa} we can obtain
\begin{align}
& \left| \vp_n^\top \left( \p_\theta\mbbh_n(\tes) - \p_\theta\ell_n(\tes)\right) \right|
\nn\\
&\lesssim 
\left| r_n^{-1} \p_\mu\mbbh_n(\tesnn_{n}) \right| + \left| r_n^{-1} \p_\lam\mbbh_n(\tesnn_{n}) \right| + 
\left| \frac{l'}{\sqrt{n}} \p_\beta \mbbh_n(\tesnn_{n}) \right| + \left| \frac{l'}{\sqrt{n}} \p_\sig \mbbh_n(\tesnn_{n}) \right|
\nn\\
&\lesssim O_p\left(\sqrt{n}\, h^{1+\kappa - 1/\beta_0} \vee \sqrt{n}\,h^\kappa\right)
+ O_p\left(\sqrt{n}\, h^{1+\kappa - 1/\beta_0} \vee \sqrt{n}\,h^\kappa\right)
\nn\\
&{}\qquad + O_p\left(\sqrt{n}\, h^{1+\kappa - 1/\beta_0} (l')^C\right) + O_p\left(\sqrt{n}\, h^{1+\kappa - 1/\beta_0} (l')^C\right)
\nn\\
&\lesssim 
O_p\big(n^{-\del_1} \vee n^{1/2-\kappa}\big) \cip 0.
\nonumber
\end{align}
%
This concludes that $R'_{0,n}=o_p(1)$.
}


\medskip

\v2{
\paragraph{\textit{Step 3.}}
Let $R''_{0,n}=:(R''_{0,a,n},R''_{0,b,n}) \in \mbbr^{q+1}\times\mbbr^2$.
The goal of this step is to show $R''_{0,a,n} = o_p(1)$ and $R''_{0,b,n} = O_p(n^{1/2-r} (l')^C)$;
at this stage, the latter component may not be stochastically bounded if $r\le 1/2$ (recall \eqref{r_set}).
We have $R''_{0,n} = A_{0,n} H_{0,n}$, where
\begin{align}
A_{0,n} &:= \vp_n^\top 
\begin{pmatrix}
\p_{a}^{2}\mbbh_n(\tesnn'_{0,n}) - \p_{a}^{2}\mbbh_n(\tesnn_{0,n}) & \text{sym.}
\\
\p_{a}\p_{b}\mbbh_n(\tesnn'_{0,n}) & \p_{b}^{2}\mbbh_n(\tesnn'_{0,n}) - \p_{b}^{2}\mbbh_n(\tesnn_{0,n})
\end{pmatrix}\vp_n ,
\nn\\
H_{0,n} &:= \vp_{n}^{-1}(\tesnn_{0,n} - \tes).
\nonumber
\end{align}
Under the assumption $\vp_{0,n}^{-1}(\tesnn_{0,n}-\tz)=O_p(1)$, recalling the block-diagonal forms \eqref{hm:vp_def} and \eqref{vp0.rate}, we see that
\begin{equation}
H_{0,n} = \vp_{n}^{-1}\vp_{0,n} \, 
\vp_{0,n}^{-1}(\tesnn_{0,n} -\tz) - \vp_{n}^{-1}(\tes -\tz)=
\begin{pmatrix}
O_p(1) \\ O_p(n^{(1-r)/2} l')
\end{pmatrix},
\label{v2-4}
\end{equation}
where the components $O_p(1)\in\mbbr^{q+1}$ and $O_p(n^{(1-r)/2} l') \in\mbbr^2$;
here and in what follows, we use the stochastic-order symbols for random variables of different dimensions, which will not cause any confusion.
}

\v2{
Turning to $A_{0,n}$, we will show that all the components of $A_{0,n}$ are at most $O_p\big( n^{-r/2} (l')^C\big)$:
\begin{equation}
A_{0,n} = O_p\big( n^{-r/2} (l')^C\big).
\label{v2-5}
\end{equation}
For the diagonal parts of $A_{0,n}$, from the same arguments as in proving \eqref{v2-6} and \eqref{v2-7} with the assumption $\vp_{0,n}^{-1}(\tesnn_{0,n}-\tz)=O_p(1)$, it holds that
\begin{align}
\left|r_n^{-2} \left(\p_a^2\mbbh_n(\tesnn'_{0,n}) - \p_a^2\mbbh_n(\tesnn_{0,n})\right)\right|
+
\left|\frac{1}{n} \tilde{\vp}_n^{\top} \left(\p_b^2\mbbh_n(\tesnn'_{0,n}) - \p_b^2\mbbh_n(\tesnn_{0,n}) \right) \tilde{\vp}_n\right|
= O_p\big( n^{-r/2} (l')^C\big).
\nonumber
\end{align}
Write $\theta=(\theta_l)_{l=1}^{p}$ and so on, and also let $\p_{a}\p_{b}\mbbh_n(\tesnn'_{0,n}) \in \mbbr^{2}\times \mbbr^{q+1}$ for the size of the matrix.
Then, for the non-diagonal part of $A_{0,n}$, we expand it as follows:
\begin{align}
\frac{1}{r_n \sqrt{n}}\tilde{\vp}_n^\top \p_{a}\p_{b}\mbbh_n(\tesnn'_{0,n}) &=: 
\frac{1}{r_n \sqrt{n}}\p_{a}\p_{b}\mbbh_n(\tz) + \sum_{l=1}^{p} \left(\frac{(h^{1-1/\beta_0})^{-1}}{n}\p_{\theta_l}\p_{a}\p_{b}\mbbh_n(\tesnn''_{0,n})\right) (\tesnn'_{0,n,l} - \theta_{0,l}).
\nonumber
\end{align}
As in the previous diagonal case, the second term on the right-hand side equals $O_p(n^{-r/2} (l')^C)$.
As for the first term, we write
\begin{equation}
\frac{1}{r_n \sqrt{n}}\p_{a}\p_{b}\mbbh_n(\tz)
= \frac{1}{r_n \sqrt{n}}\p_{a}\p_{b}\ell_n(\tz) + 
\frac{1}{r_n \sqrt{n}}\p_{a}\p_{b}\left( \mbbh_n(\tz) - \ell_n(\tz) \right).
\nonumber
\end{equation}
We have seen the explicit expressions of the components of $\p_\theta^2\ell_n(\theta)$ in Section \ref{sec_LA-proof}.
Based on them, it can be seen that all the components of $r_n^{-1}n^{-1/2}\p_{a}\p_{b}\ell_n(\tz)$ take the form
\begin{equation}
\frac1n \sumj \pi_{j-1}(\tz)\psi(\ep_j(\tz)) + O(h^2)
\nonumber
\end{equation}
for some $\mcf_{t_{j-1}}$-measurable random variable $\pi_{j-1}(\tz)$ such that $|\pi_{j-1}(\tz)|\lesssim (1+|Y_{t_{j-1}}|)(l')^C$  
and for some odd function $\psi$ (hence $\E[\psi(\ep_j(\tz))]=0$); the last term ``$O(h^2)$'' only appears in $\p_\lam\p_\beta\ell_n(\theta)$.
Burkholder's inequality for the martingale difference arrays gives $n^{-1}\sumj \pi_{j-1}(\tz)\psi(\ep_j(\tz)) = O_p(n^{-1/2}(l')^{C})$.
We conclude that $r_n^{-1}n^{-1/2}\p_{a}\p_{b}\ell_n(\tz) = O_p(n^{-1/2}(l')^{C})$.
Next, we write $\mbbh_n(\theta) - \ell_n(\theta) = \sumj B_j(\theta)d_{\ep,j}(\theta)$ for the expression \eqref{v2-3}.
The following estimates hold:
$|d_{\ep,j}(\theta)| \lesssim h^{1+\kappa-1/\beta}$,
$|\p_a \p_b d_{\ep,j}(\theta)| \lesssim h^{1+\kappa-1/\beta}(1+l')$,
$|B_j(\theta)| \lesssim 1$,
$|\p_a B_j(\theta)| \lesssim (1+|Y_{t_{j-1}}|)h^{1-1/\beta}$,
$|\p_b B_j(\theta)| \lesssim 1+l'$,
and $|\p_a\p_b B_j(\theta)| \lesssim (1+l')(1+|Y_{t_{j-1}}|)h^{1-1/\beta}$.
Therefore, by \eqref{v2-1},
\begin{align}
\left| \frac{1}{r_n \sqrt{n}}\p_{a}\p_{b}\left( \mbbh_n(\tz) - \ell_n(\tz) \right) \right|
&= \left| \frac{1}{r_n \sqrt{n}}\sumj \left.\p_{a}\p_{b} \left( B_j(\theta)d_{\ep,j}(\theta) \right)\right|_{\theta=\tz} \right|
\nn\\
&\lesssim (1+l') h^{1+\kappa-1/\beta_0} \,\frac1n \sumj (1+|Y_{t_{j-1}}|)
\nn\\
&= O_p\left(\frac{(l')^C}{\sqrt{n}}\right) \sqrt{n}\,h^{1+\kappa-1/\beta_0}
=o_p\left(\frac{(l')^C}{\sqrt{n}}\right).
\nonumber
\end{align}
Since $r\le 1$, we have concluded \eqref{v2-5}.}

\v2{
The desired stochastic orders follows from \eqref{v2-4} and \eqref{v2-5}:
\begin{align}
R''_{0,n} = A_{0,n} H_{0,n} 
= O_p\big( n^{-r/2} (l')^C\big)
\begin{pmatrix}
O_p(1) \\ O_p(n^{(1-r)/2} l')
\end{pmatrix}
=
\begin{pmatrix}
o_p(1) \\ O_p\left(n^{1/2-r} (l')^C\right))
\end{pmatrix}.
\label{v2-8}
\end{align}
}

\medskip

\v2{
\paragraph{\textit{Step 4.}}
We are now able to derive the convergence rate of $\tesnn_{1,n}-\tes$.
Recall the definition \eqref{def_K} of $K\in\mbbn$ and the initial rate of convergence \eqref{vp0.rate}.
\begin{itemize}
\item First, we consider $r>1/2$.
Then, $R''_{0,n}=o_p(1)$ from \eqref{v2-8}, so that we can take $\vp_{1,n}=\vp_n$:
by Steps 1 to 3 and \eqref{NR.proof-1}, $\vp_n^{-1}(\tesnn_{1,n}-\tes)=o_p(1)$.
This means that a single iteration is enough if we can take $r>1/2$ from the beginning.
\item Turning to $r\in(0,1/2]$, we pick a constant $\ep'\in(0,r/2)$ (hence $r-\ep'>r/2$), which is to be taken sufficiently small later.
Define
\begin{equation}
\vp_{1,n}=\vp_{1,n}(\ep') := \diag\left(
r_{n}^{-1}I_{q+1},\, n^{-(r-\ep')}
\begin{pmatrix}
1 & 0 \\ 0 & l'
\end{pmatrix}
\right).
\nn
\end{equation}
Again by Steps 1 to 3 and \eqref{NR.proof-1},
$\vp_{1,n}^{-1}\vp_n \hat{\mci}_{0,n}^{-1}=\diag(O_p(1),O_p\big((l')^C n^{r-\ep'-1/2}\big))$ and
\begin{align}
\vp_{1,n}^{-1}(\tesnn_{1,n}-\tes)
&=
\begin{pmatrix}
O_p(1) & O
\\
O & O_p\big((l')^C n^{r-\ep'-1/2}\big)
\end{pmatrix}
\left\{o_p(1) + 
\begin{pmatrix}
o_p(1) \\ O_p\left(n^{1/2-r} (l')^C\right)
\end{pmatrix}
\right\}
\nn\\
&= o_p(1) + 
\begin{pmatrix}
o_p(1) \\ O_p\big(n^{-\ep'} (l')^C\big)
\end{pmatrix}
=o_p(1).
\nonumber
\end{align}
It follows that the rate of convergence for estimating $(\beta,\sig)$ gets improved from $\diag(n^{r/2}, n^{r/2}/ l')$ of $\tesnn_{0,n}$ to $\diag(n^{r-\ep'}, n^{r-\ep'}/ l')$ of $\tesnn_{1,n}$;
this can be seen as a matrix-norming counterpart of the (near-)doubling phenomenon in the one-step estimation; see for example \cite[Section 5.5]{Zac71}.
To improve the rate further, we apply \eqref{tes.k_def} to obtain $\tesnn_{2,n}$ from $\tesnn_{1,n}$, so that the rate of convergence for estimating $(\beta,\sig)$ gets improved from $\diag(n^{r-\ep'}, n^{r-\ep'}/ l')$ to $\diag(n^{2r-3\ep'}, n^{2r-3\ep'}/ l')$; here again, we can control the constant $\ep'>0$ to be sufficiently small.
This procedure is iterated $K-1$ times, resulting in the rate $\diag(n^{2^{K-2}r-\ep'_0}, n^{2^{K-2}r-\ep'_0}/ l')$ with $\ep'_0$ being small enough to ensure that $2(2^{K-2}r-\ep'_0)>1/2$. Then, the last ($K$th-step) application of \eqref{tes.k_def} is the same as in the case of $r>1/2$ mentioned above.
\end{itemize}
These observations conclude \eqref{NR.proof-5}.
}

\medskip

Thus, we have arrived at the following claim.

\v2{
\begin{thm}
\label{thm_Kstep}
Suppose that $\tesnn_{0,n}$ satisfies that $\vp_{0,n}^{-1}(\tesnn_{0,n} - \tz) = O_p(1)$ with \eqref{vp0.rate} and \eqref{r_set}, and define $K$ as in \eqref{def_K}.
Then, the $K$-step estimator $\tesnn_{K,n}$ defined through \eqref{tes.k_def} satisfies \eqref{NR.proof-5}, hence asymptotically efficient (by Theorem \ref{thm_local.LF.asymp}):
\begin{equation}
\vp_n(\tz)^{-1}(\tesnn_{K,n} - \tz) = \mci_{n}(\tz)^{-1}\D_n(\tz) + o_{p}(1) \cil MN_{p,\tz}\left(0,\, \mci(\tz)^{-1} \right)
\label{thm_onestep-1}
\end{equation}
\end{thm}
}

\medskip

Because of the diagonality of $\vp_{0,n}$, Theorem \ref{thm_Kstep} makes it possible to construct an initial estimator $\tesnn_{0,n}=(\hat{\lam}_{0,n},\hat{\mu}_{0,n},\hat{\beta}_{0,n},\hat{\sig}_{0,n})$ individually for each component.

Having \eqref{thm_onestep-1} in hand, we can construct consistent estimators $\hat{\mci}_{\lam,\mu,n}\cip \mci_{\lam,\mu}(\tz)$ and $\hat{\mci}_{\beta,\sig,n}\cip \mci_{\beta,\sig}(\tz)$, and then prove the Studentization:
\begin{equation}
\left( \hat{\mci}_{\lam,\mu,n}^{1/2}\sqrt{n} \,h^{1-1/\besnn_{K,n}} \binom{\lesnn_{K,n}-\lam_0}{\mesnn_{K,n}-\mu_0},
~\hat{\mci}_{\beta,\sig,n}^{1/2}\sqrt{n}\,\tilde{\vp}_n(\tesnn_{K,n})^{-1}\binom{\besnn_{K,n}-\beta_0}{\sesnn_{K,n}-\sig_0}
\right) \cil N_{p}(0,I_{p}).
\label{asymp.norm}
\end{equation}
Indeed, this follows by noting the following facts.

\begin{itemize}
\item For construction of $\hat{\mci}_{\lam,\mu,n}$ and $\hat{\mci}_{\beta,\sig,n}$:
\begin{itemize}
\item In the expressions \eqref{hm:FIm_def1} and \eqref{hm:FIm_def2}, we can replace the (Riemann) $dt$-integrals by the corresponding sample quantities:
\begin{equation}
\frac1n \sumj \big( Y_{t_{j-1}}^2, Y_{t_{j-1}}X_{t_{j-1}}\big) \cip \frac1T \int_0^T \big( Y_{t}^2, Y_{t}X_{t}\big) dt.
\nn
\end{equation}

\item The elements of the form $\E_{\tz}[H(\ep;\beta_0)] = \int H(\ep;\beta_0)\phi_{\beta_0}(\ep)d\ep$ with $H(\ep;\beta)$ smooth in $\beta$ can be evaluated through a numerical integration involving the density $\phi_\beta(\ep)$ and its partial derivatives with respect to $(\beta,\ep)$, with plugging-in the estimate $\besnn_{K,n}$ for the value of $\beta$ (the initial estimator $\besnn_{0,n}$ is enough).

\item Again note that $n^{v}(\besnn_{K,n} - \beta_0,\, \sesnn_{K,n} - \sig_0) = o_p(1)$ for any sufficiently small $v \in (0,1/2)$, so that $h^{1-1/\besnn_{K,n}} / h^{1-1/\beta_{0}} = (1/h)^{1/\besnn_{K,n}-1/\beta_{0}} \cip 1$.
The values $\overline{\vp}_{lm}(\tz)$ contained in $\mci_{\beta,\sig}(\tz)$ are estimated by plugging-in $\tesnn_{K,n}$ in \eqref{hm:vp-conditions}:
\begin{equation}
\left\{\begin{array}{l}
\besnn_{K,n}^{-2} l' \vp_{11,n}(\tesnn_{K,n}) + \sesnn_{K,n}^{-1}\vp_{21,n}(\tesnn_{K,n}) \cip \overline{\vp}_{21}(\tz), \nn\\
\besnn_{K,n}^{-2} l' \vp_{12,n}(\tesnn_{K,n}) + \sesnn_{K,n}^{-1}\vp_{22,n}(\tesnn_{K,n}) \cip \overline{\vp}_{22}(\tz), \nn\\
\vp_{11,n}(\tesnn_{K,n}) \cip \overline{\vp}_{11}(\tz), \nn\\
\vp_{12,n}(\tesnn_{K,n}) \cip \overline{\vp}_{12}(\tz). \nn\\
\end{array}\right.
\nonumber
\end{equation}
We can replace $(\besnn_{K,n},\sesnn_{K,n})$ by $(\besnn_{0,n},\sesnn_{0,n})$ all through the above.

\end{itemize}

\item Since $\vp_{n}^{-1}(\tesnn_{K,n}-\tz)=O_p(1)$, it follows that
\begin{align}
\sqrt{n}\,\tilde{\vp}_n(\tesnn_{K,n})^{-1}\binom{\besnn_{K,n}-\beta_0}{\sesnn_{K,n}-\sig_0}
&= \sqrt{n}\, \left(\tilde{\vp}_n(\tz)^{-1} + O_p\big( (l')^C n^{-1/2}\big)\right)
\binom{\besnn_{K,n}-\beta_0}{\sesnn_{K,n}-\sig_0}
\nn\\
&= \sqrt{n}\, \tilde{\vp}_n(\tz)^{-1} \binom{\besnn_{K,n}-\beta_0}{\sesnn_{K,n}-\sig_0} + O_p\big( (l')^C n^{-1/2}\big)
\nn\\
&= \sqrt{n}\, \tilde{\vp}_n(\tz)^{-1} \binom{\besnn_{K,n}-\beta_0}{\sesnn_{K,n}-\sig_0} + o_p(1).
\nonumber
\end{align}
\end{itemize}
The property \eqref{asymp.norm} entails
\begin{equation}
\left| \hat{\mci}_{\lam,\mu,n}^{1/2}\sqrt{n} \,h^{1-1/\besnn_{K,n}} \binom{\lesnn_{K,n}-\lam_0}{\mesnn_{K,n}-\mu_0} \right|^2
+
\left| \hat{\mci}_{\beta,\sig,n}^{1/2}\sqrt{n}\,\tilde{\vp}_n(\tesnn_{K,n})^{-1}\binom{\besnn_{K,n}-\beta_0}{\sesnn_{K,n}-\sig_0} \right|^2
\cil \chi^2(p)=\chi^2(q+3),
\nonumber
\end{equation}
which can be used for constructing an approximate confidence ellipsoid and for goodness-of-fit testing, in particular variable selection among the components of $X$.

\medskip

\begin{rem}\normalfont
From the proof of Theorem \ref{thm_Kstep}, we see that it is possible to weaken \eqref{beta&kappa} as $\beta_0 > 2/3$ if the integrated-process sequence $(\int_j X_s ds)_{j=1}^n$ is observable.
\v2{
Moreover, It is possible to remove \eqref{beta&kappa} if the model is the Markovian $Y_t = Y_0 + \int_0^t (\mu-\lam Y_s)ds + \sig J_t$ with constant $\mu\in\mbbr$ with modifying the definition \eqref{epj'_def} as in the estimating function of \cite{CleGlo20}.
However, we worked under \eqref{X-kappa} and \eqref{epj'_def} to deal with a possibly time-varying $X$.
}
\end{rem}

\v2{
\begin{rem}\normalfont
\label{rem_block.diag.form}
The standard form of the one-step estimator is not \eqref{tes.k_def}, but
\begin{equation}
\tesnn_{k,n} = \tesnn_{k-1,n} + \left( -\p_\theta^2\mbbh_n(\tesnn_{k-1,n})\right)^{-1} \p_\theta\mbbh_n(\tesnn_{k-1,n}).
\nonumber
\end{equation}
By inspecting the proof of Theorem \ref{thm_Kstep}, we found that the off-block-diagonal part $-\p_a \p_b \mbbh_n(\tesnn_{k-1,n})$ made the claim therein invalid.
This has happened since the rate of convergence for estimating the component $b=(\beta,\sig)$ could be too slow.
Still, because of the block-diagonality of the original form \eqref{hm:FIm_def0}, it seems to be a natural and reasonable strategy to use the block-diagonal form from the beginning of defining \eqref{tes.k_def}.
\end{rem}
}

\begin{rem}\normalfont
\label{rem:rates.PLDI}
The necessity of more than one iteration ($K\ge 2$) would be a technical one.
If we could verify the tail-probability estimate $\sup_n \pr[|r_n(\lesnn_{0,n}-\lam_0,\mesnn_{0,n}-\mu_0)| \ge s]\lesssim s^{-M}$ for a sufficiently \textit{large} $M>0$, 
\v2{
then it is possible to deduce the optimality of the one-step Newton-Raphson procedure even when a strategy of construction $(\besnn_{0,n},\sesnn_{0,n})$ is not smooth in $(\lesnn_{0,n},\mesnn_{0,n})$ as in the function $\hat{M}_n(a')$ in Section \ref{sec:mm}.
}
However, the model under consideration is heavy-tailed and it seems impossible to deduce such a bound since we cannot make use of the localization for that purpose.
\end{rem}


\subsection{Specific preliminary estimators}

In this section, we consider a specific construction of $\tesnn_{0,n}=(\lesnn_{0,n},\mesnn_{0,n},\besnn_{0,n},\sesnn_{0,n})$ satisfying $\vp_{0,n}^{-1}(\tesnn_{0,n} - \tz) = O_p(1)$ with $\vp_{0,n}$ given by \eqref{vp0.rate}.
We keep assuming that available sample is $\{(X_{t_j}, Y_{t_j})\}_{j=0}^{n}$ and the conditions \eqref{X-kappa} and \eqref{beta&kappa} are in force.
We will proceed in two steps.
\begin{enumerate}
\item First, we will estimate the trend parameter $(\lam,\mu)$ by the least absolute deviation (LAD) estimator, which will turn out to be rate-optimal, and asymptotically mixed-normally distributed; although the identification of the asymptotic distribution is not necessary here, it would be of independent interest (see Section \ref{appendix_LAD_AMN}).
\item Next, by plugging in the LAD estimator we construct a sequence of residuals for the noise term, based on which we will consider the lower-order fractional moment matching.
\end{enumerate}
Recall that we are working under the localization \eqref{localization.moment} by removing large jumps of $J$.

\subsubsection{LAD estimator}
\label{sec:LAD}

Let us recall the autoregressive structure together \eqref{Y_sol} with the approximation of the (non-random) integral:
\begin{align}
Y_{t_j} & = e^{-\lam_0 h}Y_{t_{j-1}} + \mu_0 \cdot \zeta_{j}(\lam_0)h + \sig_0 \int_j e^{-\lam_0 (t_j -s)}dJ_s
\nn\\
&= Y_{t_{j-1}} - \lam_0 h Y_{t_{j-1}}  + \mu_0 \cdot X_{t_{j-1}} h + \sig_0 \int_j e^{-\lam_0 (t_j -s)}dJ_s + h^{1/\beta_0} \del'_{j-1},
\label{LAD-1}
\end{align}
where
\begin{equation}
\del'_{j-1}=\del'_{j-1}(\tz) := h^{-1/\beta_0}  \left( Y_{t_{j-1}} (
\v2{e^{-\lam_0 h}} - 1 + \v2{\lam_0} h) + \v2{\mu_0} \cdot \left(\zeta_{j}(\v2{\lam_0}) - X_{t_{j-1}}\right) h \right)
\nonumber
\end{equation}
 is an $\mcf_{t_{j-1}}$-measurable random variable such that
 \begin{equation}
|\del'_{j-1}| \lesssim (1+|Y_{t_{j-1}}|) h^{1+\kappa -1/\beta_0}.
\label{LAD-5}
\end{equation}

We define the LAD estimator $(\lesnn_{0,n},\mesnn_{0,n}) \in \mbbr^{q+1}$ by any element $(\lesnn_{0,n},\mesnn_{0,n}) \in \argmin_{(\lam,\mu)} M_n(\lam,\mu)$ with leaving $(\beta,\sig)$ unknown, where
\begin{align}
M_n(\lam,\mu) &:= \sumj \left| Y_{t_j}- Y_{t_{j-1}} - \left( - \lam Y_{t_{j-1}} + \mu \cdot X_{t_{j-1}} \right) h \right|.
\label{LAD_obj.func}
\end{align}
This is a slight modification of the previously studied approximate LAD estimator in \cite{Mas10ejs} concerning the ergodic locally stable OU process.

We introduce the following convex random function on $\mbbr\times \mbbr^q$ (recall the notation \eqref{def:r_n}):
\begin{equation}
\Lam_{n}(u,v) := 
\frac{1}{\sig_0 \eta(\lam_0\beta_0 h)^{1/\beta_0}\, h^{1/\beta_0}}
\left\{M_n\left( \lam_0 + \frac{u}{r_n},\, \mu_0 + \frac{1}{r_n}v\right) - M_n(\lam_0, \mu_0) \right\}.
\nonumber
\end{equation}
The minimizer of $\Lam_n$ is $\hat{w}_n:=(\hat{u}_n,\hat{v}_n)$ where $\hat{u}_n := r_n(\lesnn_{0,n}-\lam_0)$ and $\hat{v}_n := r_n(\mesnn_{0,n}-\mu_0)$.
Further, letting $z_{j-1}:=(-Y_{t_{j-1}}, X_{t_{j-1}})$, $w:=(u,v)$, and
\begin{equation}
\ep'_j := \frac{1}{\eta(\lam_0\beta_0 h)^{1/\beta_0}\, h^{1/\beta_0}} \int_j e^{-\lam_0 (t_j -s)}dJ_s
~\stackrel{\pr_{\tz}}{\sim}~\text{i.i.d.}~\mcl(J_1),
\nonumber
\end{equation}
we also introduce the quadratic random function
\begin{equation}
\Lam_{n}^\sharp(w) := \D'_{n}[w] + \frac{1}{2}\Gam_{0}[w,w],
\nonumber
\end{equation}
where
\begin{align}
\D'_{n} &:= \v2{-} \sumj \frac{1}{s_{0,n}\sqrt{n}} \, \sgn\left(\ep'_j + \del_{j-1}'\right) z_{j-1}, \nn\\
\nonumber \\
\Gam_0 &:= \frac{2\phi_{\beta_0}(0)}{\sig_0^2}\frac1T\int_0^T
\begin{pmatrix}
Y_{t}^2 & -Y_{t}X_{t}^\top \\ -Y_{t}X_{t} & X_{t}^{\otimes 2}
\end{pmatrix}dt,
\nonumber
\end{align}
where $s_{0,n}:=\sig_0 \eta(\lam_0\beta_0 h)^{1/\beta_0} = (1+o(1)) \sig_0$.
The a.s. positive definiteness of $\Gam_{0}$ (see Section \ref{sec:likelihood.asymptotics}) implies that $\argmin\Lam_{n}^{\sharp}$ a.s. consists of the single point $\hat{w}_{n}^{\sharp}:=-\Gam_{0}^{-1}\D'_{n}$.
Then, our objective is to prove that
\begin{equation}
\hat{w}_{n}=\hat{w}_{n}^{\sharp}+o_{p}(1).
\label{LAD-2}
\end{equation}
The proof is analogous to \cite[Proof of Theorem 2.1]{Mas10ejs}, hence we will appropriately omit the full technical details, referring to the corresponding parts therein.

\medskip

By \eqref{LAD-1} and \eqref{LAD_obj.func}, we have
\begin{align}
\Lam_{n}(w)
&= \sumj \left(\left| \ep'_j + \v2{\frac{\del_{j-1}'}{s_{0,n}}} 
- \frac{1}{s_{0,n}\sqrt{n}} w\cdot z_{j-1}\right| - \left| \ep'_j + \v2{\frac{\del_{j-1}'}{s_{0,n}}}\right| \right).
\nonumber
\end{align}
As in \cite[Eq.(4.6)]{Mas10ejs}, we can write $\Lam_n(w) = \D'_n[w] + Q_n(w)$ where
\begin{align}
Q_{n}(w) &:= 2\sumj \int_{0}^{w\cdot z_{j-1} / (s_{0,n}\sqrt{n})}
\left\{I\left(\ep'_j + \v2{\frac{\del_{j-1}'}{s_{0,n}}} \le s\right) - I\left(\ep'_j + \v2{\frac{\del_{j-1}'}{s_{0,n}}} \le 0\right)\right\}ds.
\nonumber
\end{align}
For a moment, let us suppose that
\begin{align}
\D'_n &= O_p(1),
\label{LAD-3.1} \\
Q_n(w) &= \frac12 \Gam_0 [w,w] + o_p(1), \qquad w\in\mbbr^{1+q}.
\label{LAD-3.2}
\end{align}
Then, we can make use of the argument of \cite{HjoPol93} to conclude \eqref{LAD-2}.
To see this, we note the inequality due to \cite[Lemma 2]{HjoPol93}:
for any $\ep>0$,
\begin{equation}
\pr\left[ |\hat{w}_{n}-\hat{w}_{n}^{\sharp}|\ge\ep \right]
\le 
\pr\left[ \sup_{w:\, |w-\hat{v}_{n}^{\sharp}|\le\ep}|\del_{n}(\v2{w})| 
\ge 
\frac{1}{2} \left(\inf_{(w,z):\, |z|=1, \atop w=\hat{w}_{n}^{\sharp}+\ep z}\Lam_{n}^{\sharp}(\v2{w}) - \Lam_{n}^{\sharp}(\hat{w}_n^\sharp)\right) \right],
\nonumber
\end{equation}
where $\del_n(w) := \Lam_n(w) - \Lam_n^\sharp(w)$.
Obviously, $\Lam_{n}^{\sharp}(\hat{w}_n^\sharp) = \v2{-} (1/2)\D'_n \cdot \Gam_0^{-1}\D'_n$.
By straightforward computations, we obtain
\begin{align}
\inf_{(w,z):\, |z|=1, \atop w=\hat{w}_{n}^{\sharp}+\ep z}\Lam_{n}^{\sharp}(\v2{w}) - \Lam_{n}^{\sharp}(\hat{w}_n^\sharp)
\ge \ep^2 \lam_{\min}(\Gam_0).
\nonumber
\end{align}
Also, because of the convexity, we have the uniform convergence $\sup_{w\in A}\left|\del_n(w)\right| \cip 0$ for each compact $A\subset \mbbr^{1+q}$ (see \cite[Lemma 1]{HjoPol93}).
Note that $\hat{w}_{n}^{\sharp} = O_p(1)$ by \eqref{LAD-3.1} and the a.s. positive definiteness of $\Gam_0$.
Given any $\ep,\ep'>0$, we can find sufficiently large $K>0$ and $N\in\mbbn$ for which the following three estimates hold simultaneously:
\begin{align}
& \sup_n \pr[|\hat{w}_{n}^{\sharp}| > K] < \ep'/3, \nn\\
& \sup_{n\ge N} \pr\left[
\sup_{w:\, |w|\le K + \ep}|\del_{n}(v)| > \ep' \right] < \ep'/3, \nn\\
& \pr\left[ \ep' \ge \frac{\ep^2}{2} \lam_{\min}(\Gam_0)\right] < \ep'/3.
\nonumber
\end{align}
Piecing together the above arguments concludes that, for any $\ep,\ep'>0$, there exists an $N\in\mbbn$ such that
$\sup_{n\ge N} \pr\left[ |\hat{w}_{n}-\hat{w}_{n}^{\sharp}|\ge\ep \right] < \ep'$.
This establishes \eqref{LAD-2}, and it follows that
\begin{equation}
\hat{w}_{n} = -\Gam_{0}^{-1}\D'_{n} + o_p(1) = O_p(1).
\label{LAD-4}
\end{equation}
It remains to prove \eqref{LAD-3.1} and \eqref{LAD-3.2}.
Below, we will write $\pr^{j-1}$ and $\E^{j-1}$ for the conditional probability and expectation given $\mcf_{t_{j-1}}$, respectively.

\medskip

\noindent
\textit{Proof of \eqref{LAD-3.1}}
follows on showing $\D_n(1)=O_p(1)$ and $R_{1,n}=o_p(1)$, where
\begin{align}
\D_{n}(t) &:= \sum_{j=1}^{[nt]} \frac{1}{\sig_{0}\sqrt{n}} \, 
\left\{ \sgn\left(\ep'_j + \v2{\frac{\del_{j-1}'}{s_{0,n}}}\right) - \E_{0}^{j-1}\left[ \sgn\left(\ep'_j + \v2{\frac{\del_{j-1}'}{s_{0,n}}}\right) \right]\right\} z_{j-1}, \qquad t\in[0,1], \nn\\
R_{1,n} &:= \sumj \frac{1}{\sig_{0}\sqrt{n}} \, \E_{0}^{j-1}\left[ \sgn\left(\ep'_j + \v2{\frac{\del_{j-1}'}{s_{0,n}}}\right) z_{j-1}\right].
\label{LAD-7}
\end{align}
The (matrix-valued) predictable quadratic variation process of $\{\D_{n}(\cdot)\}_{t\in[0,1]}$ is given by
\begin{equation}
\la \D_n(\cdot)\ra_t 
:= \sig_0^{-2} \frac1n \sum_{j=1}^{[nt]} 
\left\{ \sgn\left(\ep'_j + \v2{\frac{\del_{j-1}'}{s_{0,n}}}\right) - \E_{0}^{j-1}\left[ \sgn\left(\ep'_j + \v2{\frac{\del_{j-1}'}{s_{0,n}}}\right) \right]\right\}^2
z_{j-1}^{\otimes 2}
\nonumber
\end{equation}
We apply the Lenglart inequality \cite[I.3.31]{JacShi03} for the submartingale $|\D_n(t)|^2$:
for any $K,L>0$,
\begin{align}
\sup_{n} \pr\left[ \sup_{t\in[0,1]}|\D_{n}(t)| \ge K \right]
&\lesssim \frac{L}{K} + \sup_{n} \pr\left[ \frac{1}{n}\sumj |z_{j-1}|^2 \ge C \sig_0^2 L\right]
\nn\\
&\lesssim \frac{L}{K} + \sup_{n} \pr\left[ \frac{1}{n}\sumj (1+|Y_{t_{j-1}}|)^2 \ge C L\right].
\nn
\end{align}
We have $n^{-1}\sumj (1+|Y_{t_{j-1}}|)^2=O_p(1)$. To conclude that $\D_n:=\D_n(1)=O_p(1)$, let $L$ and $K$ sufficiently large in this order.
To see $R_{1,n}=o_p(1)$, we proceed in exactly the same way as in \cite[pp.544--545]{Mas10ejs}:
by partly using \eqref{beta&kappa} and \eqref{LAD-5},
\begin{align}
|R_{1,n}| &= \left| \frac1n \sumj 2\sqrt{n} \, z_{j-1} \int_0^{\v2{\del'_{j-1}/s_{0,n}}} \phi_{\beta_0}(y)dy \right| \nn\\
&\lesssim \frac1n \sumj \sqrt{n} |z_{j-1}| |\del'_{j-1}| 
\lesssim \frac1n \sumj (1+|Y_{t_{j-1}}|)^2 \sqrt{n} h^{1+\kappa-1/\beta_0} \nn\\
&=O_p\left(\sqrt{n} h^{1+\kappa-1/\beta_0}\right) = O_p\left(h^{1/2+\kappa-1/\beta_0}\right) = o_p(1).
\nonumber
\end{align}
Thus we have obtained \eqref{LAD-3.1}, and now we can replace $\D_n'$ by $\D_n$ in \eqref{LAD-4}:
\begin{equation}
\hat{w}_{n} = -\Gam_{0}^{-1}\D_{n} + o_p(1) = O_p(1).
\label{LAD-6}
\end{equation}

\medskip

\noindent
\textit{Proof of \eqref{LAD-3.2}.}
We decompose $Q_{n}(w) =: \sumj \zeta_j(w)$ as $Q_n(w)=Q_{1,n}(w) + Q_{2,n}(w)$, where 
$Q_{1,n}(w) := \sumj \E^{j-1}[\zeta_j(w)]$ and $Q_{2,n}(w) := \sumj (\zeta_j(w) - \E^{j-1}[\zeta_j(w)])$.
Then, for each $w\in\mbbr^{1+q}$ we can readily mimic the flow of \cite[pp.545--546]{Mas10ejs} (for handling the term $\mathbb{Q}_n(u)$ therein).
The sketches are given below.
\begin{itemize}
\item We have
\begin{equation}
Q_{1,n}(w) = \frac12 \Gam_n [w,w] + A_n(w),
\nonumber
\end{equation}
where
\begin{equation}
\Gam_n := 
\frac{2\phi_{\beta_0}(0)}{s_{0,n}^2}\frac{1}{n}\sum_{j=1}^{n}
\begin{pmatrix}
Y_{t_{j-1}}^2 & -Y_{t_{j-1}}X_{t_{j-1}}^\top \\ -Y_{t_{j-1}}X_{t_{j-1}} & X_{t_{j-1}}^{\otimes 2}
\end{pmatrix}
= \Gam_0 + o_p(1),
\nonumber
\end{equation}
and where
\begin{align}
|A_n(w)| &\lesssim 
\left|\frac{1}{n}\sumj (w\cdot z_{j-1})^{2}\left\{
\phi_{\beta_0}\left(-\v2{\frac{\del_{j-1}'}{s_{0,n}}}\right) - \phi_{\beta_0}(0)\right\}\right| \nonumber \\
&{}\qquad + \left| \sumj \int_{0}^{w\cdot z_{j-1}/(s_{0,n}\sqrt{n})}s^{2}\int_{0}^{1}(1-y)
\p \phi_{\beta_0}\left(sy-\v2{\frac{\del_{j-1}'}{s_{0,n}}}\right)dyds \right|
\nn\\
&\lesssim \frac1n \sumj (1+|Y_{t_{j-1}}|)^4 (1+|w|)^4
\left( h^{2(1+\kappa-1/\beta_0)} \vee \frac1n \vee \frac{h^{1+\kappa-1/\beta_0}}{\sqrt{n}} \right)
=o_p(1).
\nonumber
\end{align}

\item We have $Q_{2,n}(w)=o_p(1)$: by the Burkholder-Davis-Gundy inequality,
\begin{align}
& \E\left[\left( \sumj (\zeta_j(w) - \E^{j-1}[\zeta_j(w)]) \right)^2 \right] \nn\\
&\lesssim \sumj \E\left[\left(
\int_{0}^{|w\cdot z_{j-1} / (s_{0,n}\sqrt{n})|} I\left(\left|\ep'_j + \v2{\frac{\del_{j-1}'}{s_{0,n}}}\right| \le s\right) ds
\right)^2 \right] \nn\\
&\lesssim \sumj \frac{|w|}{\sqrt{n}}
\E\left[ |z_{j-1}| \int_{0}^{|w\cdot z_{j-1} / (s_{0,n}\sqrt{n})|} \pr^{j-1}\left[ \left|\ep'_j + \v2{\frac{\del_{j-1}'}{s_{0,n}}}\right| \le s\right]ds \right]
\nn\\
&\lesssim \sumj \frac{|w|}{\sqrt{n}}\E\left[ |z_{j-1}| \int_{0}^{|w\cdot z_{j-1} / (s_{0,n}\sqrt{n})|} 
\left(s+\left|\v2{\frac{\del_{j-1}'}{s_{0,n}}}\right|\right) ds \right]
\nn\\
&\lesssim (1+|w|)^3 \frac1n \sumj \E\left[ (1+|Y_{t_{j-1}}|)^3\right] 
\left( \frac{1}{\sqrt{n}} \vee h^{1+\kappa-1/\beta_0} \right) \nn\\
&=O\left( \frac{1}{\sqrt{n}} \vee h^{1+\kappa-1/\beta_0} \right) =o(1).
\nonumber
\end{align}
\end{itemize}
Summarizing the above yields \eqref{LAD-3.2}.

The tightness \eqref{LAD-6} is sufficient for our purpose.
As a matter of fact, the LAD estimator $(\lesnn_{0,n},\mesnn_{0,n})$ is asymptotically mixed-normally distributed. We give the details in Section \ref{appendix_LAD_AMN}.

\subsubsection{Rates of convergence at the moment matching for $(\beta,\sig)$}
\label{sec:mm}

The remaining task is to construct a specific estimator $(\besnn_{0,n},\sesnn_{0,n})$ such that
\begin{equation}
\left(
n^{r/2}(\besnn_{0,n}-\beta_0),\,\frac{n^{r/2}}{l'}(\sesnn_{0,n}-\sig_0) \right) = O_{p}(1).
\label{ME-1}
\end{equation}
This can be achieved simply by fitting some appropriate moments;
for this purpose, the localization does not make sense, since precise expressions of truly existing moments without the localization come into play.
Here we consider, as in \cite{BroMas18}, the pair of the absolute moments of order $r$ and $2r$.

Let $a'\in(0, \beta_0/2)$ and define
\begin{equation}
\hat{M}_n(a') := \frac1n \sumj \left|
Y_{t_{j}} - Y_{t_{j-1}} + \lesnn_{0,n} Y_{t_{j-1}}h  - \mesnn_{0,n} \cdot X_{t_{j-1}} h
\right|^{a'}.
\nonumber
\end{equation}
Let
\begin{equation}
\ep''_j := \frac{1}{h^{1/\beta_0}} \int_j e^{-\lam_0 (t_j -s)}dJ_s
=(1+o(1))\ep'_j,
\nonumber
\end{equation}
which are approximately i.i.d. with common distribution $\mcl(J_1)$, and also let
\begin{equation}
M_n(a') := \sig_0^{a'} h^{a'/\beta_0}\frac1n \sumj \left| \ep''_j \right|^{a'}.
\nonumber
\end{equation}
We can apply the central limit theorem to ensure that 
$\sqrt{n}\left( h^{-a'/\beta_0} \sig_0^{-a'}M_n(a') - m(a';\beta_0) \right) = O_p(1)$ as soon as $a'<\beta_0/2$, where
\begin{equation}
m(a';\beta_0) := \E\big[|J_1|^{a'}\big] = \frac{2^{a'}}{\sqrt{\pi}} \frac{\Gam((a'+1)/2) \Gam(1-a'/\beta_0)}{\Gam(1-a'/2)}.
\nonumber
\end{equation}
Moreover, it follows from the discussions in Section \ref{sec:LAD} that
\begin{align}
h^{-a'/\beta_0}\sig_0^{-a'}\hat{M}_n(a')
&= \frac1n \sumj \left|
\ep''_j + \frac{1}{\sqrt{n}}\left( \sqrt{n}\del'_{j-1} - \hat{w}_n \cdot z_{j-1} \right)
\right|^{a'},
\nonumber
\end{align}
which in turn gives
\begin{align}
n^{a'/2}\left| h^{-a'/\beta_0}\sig_0^{-a'}\left(\hat{M}_n(a') - M_n(a')\right) \right| 
&\le \frac1n \sumj \left( \sqrt{n}|\del'_{j-1}| + |\hat{w}_n| |z_{j-1}| \right)^{a'} \nn\\
&\lesssim O_p\left(\sqrt{n} \, h^{1+\kappa-1/\beta_0}\right) + O_p(1) = O_p(1).
\nonumber
\end{align}
It follows that
\begin{align}
n^{a'/2} \left(h^{-a'/\beta_0} \sig_0^{-a'} \hat{M}_n(a') - m(a';\beta_0)\right) 
= O_p(1 \vee n^{(a'-1)/2}) = O_p(1).
\nonumber
\end{align}
\v2{
Now we want to take $a'=r,2r$, which necessitates that $r\in(0, \beta_0/4)$ in the current argument.
}
Then, we conclude that
\begin{equation}
n^{r/2} \left( h^{-r/\beta_0}\sig_0^{-r} \hat{M}_n(r) - m(r;\beta_0), \, 
h^{-2r/\beta_0}\sig_0^{-2r} \hat{M}_n(2r;\beta_0) - m(2r;\beta_0) \right) = O_p(1),
\nonumber
\end{equation}
so that
\begin{equation}
n^{r/2} \left( \frac{\hat{M}_n(r)^2}{\hat{M}_n(2r)} - \frac{m(r;\beta_0)^2}{m(2r;\beta_0)} \right) = O_p(1).
\nonumber
\end{equation}
There exists a bijection $f_r$ such that $f_r(m(r;\beta)^2/m(2r;\beta))=\beta$; see \cite[Section 3.2]{BroMas18} and the references therein for the related details.
Therefore, taking $\besnn_{0,n} := f_r(\hat{M}_n(r)^2 /\hat{M}_n(2r))$ results in $n^{r/2}(\besnn_{0,n} -\beta_0)=O_p(1)$, as was to be shown.
The bisection method is sufficient for finding $\besnn_{0,n}$ numerically.

Turning to $\sesnn_{0,n}$, we note that
\begin{equation}
n^{r/2} \left( \frac{h^{-r/\beta_0}\hat{M}_n(r)}{m(r;\beta_0)} - \sig_0^r \right) = O_p(1).
\label{ME-2}
\end{equation}
Let $\sesnn_{0,n} := \left(\frac{h^{-r/\besnn_{0,n}}\hat{M}_n(r)}{m(r;\besnn_{0,n})}\right)^{1/r}$:
we claim that $\frac{n^{r/2}}{l'}(\sesnn_{0,n}-\sig_0) = O_{p}(1)$.
Since $\frac{m(r;\besnn_{0,n})}{m(r;\beta_0)} = O_p(1)$,
\begin{equation}
\left| h^{r(1/\besnn_{0,n} - 1/\beta_0)} \frac{m(r;\besnn_{0,n})}{m(r;\beta_0)} -1 \right|
\le \left| h^{r(1/\besnn_{0,n} - 1/\beta_0)} - 1 \right| O_p(1) 
+ \left| \frac{m(r;\besnn_{0,n})}{m(r;\beta_0)} -1 \right|.
\nonumber
\end{equation}
Recall that $n^{r/2}(\besnn_{0,n} -\beta_0)=O_p(1)$, hence the second term in the upper bound equals $O_p(n^{-r/2})$.
As for the first term, using that $(1/\besnn_{0,n} - 1/\beta_0) l'  = O_p( l' /n^{r/2}) = o_p(1)$, we observe
\begin{align}
h^{r(1/\besnn_{0,n} - 1/\beta_0)} - 1
= \exp\left( r(1/\besnn_{0,n} - 1/\beta_0) l' \right) - 1
= O_p\left( \frac{l'}{n^{r/2}} \right) =o_p(1).
\nonumber
\end{align}
These estimates combined with \eqref{ME-2} conclude the claim: we have \eqref{ME-1} for the above constructed $(\besnn_{0,n},\sesnn_{0,n})$.
\v2{
Given an $r\in(0,\beta_0/4)$, by Theorem 3.1, the $K$-step estimator for $K> \log_2(1/r)$ is asymptotically efficient;
if $\beta_0> 1$ is supposed beforehand, then we can take an $r>1/4$ small enough to ensure that $K=2$ is enough.
}

\subsubsection{Asymptotic mixed normality of the LAD estimator}
\label{appendix_LAD_AMN}

Recall \eqref{LAD-6}: $\hat{w}_{n} = -\Gam_{0}^{-1}\D_{n} + o_{p}(1)$.
To deduce the asymptotic mixed normality, it suffices to identify the appropriate asymptotic distribution of $(\D_{n},\Gam_{0})$, equivalently of $(\D_{n},\Gam_{n})$.

First, we clarify the leading term of $\D_n$ in a simpler form.
We have $\E[\mathrm{sgn}(\ep'_{j})]=0$ and $\E[\mathrm{sgn}(\ep'_{j})^{2}]=1$, 
Observe that $\D_n = \D_{0,n} + R_{1,n} + R_{2,n}$, where $R_{1,n}$ is given in \eqref{LAD-7} and
\begin{align}
\D_{0,n} &:= \sum_{j=1}^{n} \frac{1}{\sig_{0}\sqrt{n}} \, \sgn(\ep'_j) z_{j-1}, \nn\\
R_{2,n} &:= \sum_{j=1}^{n} \frac{1}{\sig_{0}\sqrt{n}} \, 
\left( \sgn\left(\ep'_j + \del_{j-1}'\right) - \sgn(\ep'_j) \right) z_{j-1}.
\nonumber
\end{align}
We have already seen that $R_{1,n}=o_p(1)$.
We claim that $R_{2,n}=o_p(1)$. Write $R_{2,n} = \sum_{j=1}^{n} \xi_{j}$. The claim follows on showing that both $\sum_{j=1}^{n} \E^{j-1}[\xi_{j}]=o_p(1)$ and $|\sum_{j=1}^{n} \E^{j-1}[\xi_{j}^{\otimes 2}]|=o_p(1)$, but the first one obviously follows from $R_{1,n}=o_p(1)$.
The second one can be shown as follows: first, we have
\begin{align}
\sum_{j=1}^{n} \E^{j-1}[\xi_{j}^{\otimes 2}]
&= \sig_0^{-2} \frac1n \sumj \E^{j-1}\left[\left( \sgn\left(\ep'_j + \del_{j-1}'\right) - \sgn(\ep'_j) \right)^2\right] z_{j-1}^{\otimes 2}
\nn\\
&= 2\sig_0^{-2} \frac1n \sumj \left(
1 - \E^{j-1}\left[\sgn\left(\ep'_j + \del_{j-1}'\right) \sgn(\ep'_j) \right] 
\right) z_{j-1}^{\otimes 2}.
\nonumber
\end{align}
Moreover,
\begin{align}
& \E^{j-1}\left[\sgn\left(\ep'_j + \del_{j-1}'\right) \sgn(\ep'_j) \right] \nn\\
&= \left(\int_{0\vee (-\del'_{j-1})}^\infty + \int_{-\infty}^{0\vee (-\del'_{j-1})} - \int_{-\del'_{j-1}}^0 - \int_0^{-\del'_{j-1}}\right)
\phi_{\beta_0}(y)dy = 1+ D_{j-1}
\nonumber
\end{align}
for some $\mcf_{t_{j-1}}$-measurable term $D_{j-1}$ satisfying the estimate $|D_{j-1}|\lesssim |\del'_{j-1}|\lesssim (1+|Y_{t_{j-1}}|) h^{1+\kappa -1/\beta_0}$.
These observations conclude that $|\sum_{j=1}^{n} \E^{j-1}[\xi_{j}^{\otimes 2}]|=o_p(1)$.

\medskip

It remains to look at $\D_{0,n}$.
The mere convergence in distribution is not suitable for the purpose since the matrix $\Gam_0$ is random.
We will apply the weak limit theorem for stochastic integrals: we refer the reader to \cite[VI.6]{JacShi03} for a detailed account of the limit theorems as well as the standard notation used below.

We introduce the partial sum process
\begin{equation}
S^{n}_{t}:=\sum_{j=1}^{[nt]}\frac{1}{\sqrt{n}}\mathrm{sgn}(\ep'_{j}),\qquad t\in[0,1].
\nonumber
\end{equation}
We apply \cite[Lemma 4.3]{Jac07} to derive $S^{n} \scl w'$ in $\mcd(\mbbr)$ (the Skorokhod space of $\mbbr$-valued functions, equipped with the Skorokhod topology), where $w'=(w')_{t\in[0,1]}$ denotes a standard Wiener process defined on an extended probability space and independent of $\mcf$. 
Here the symbol $\scl$ denotes the ($\mcf$-)stable convergence in law, which is strictly stronger than the mere weak convergence and in particular implies the joint weak convergence in $\mcd(\mbbr^{q+2})$:
\begin{equation}
(S^{n},H^{n})\cil(w',H^{\infty})
\label{LAD.AMN_add1}
\end{equation}
for any $\mbbr^{q+1}$-valued $\mcf$-measurable {\cadlag} processes $H^{n}$ and $H^{\infty}$ such that $H^{n} \cip H^{\infty}$ in $\mcd(\mbbr^{q+1})$.

We note the following two points.

\begin{itemize}

\item We have $S^{n}\cil w'$ in $\mcd(\mbbr)$, and for each $n\in\mbbn$ the process 
$(S^{n}_{t})_{t\in[0,1]}$ is an $(\mcf_{[nt]/n})$-martingale such that $\sup_{n,t}|\D S^{n}_{t}|\le 1$. 
These facts combined with \cite[VI.6.29]{JacShi03} imply that the sequence $(S^{n})$ is predictably uniformly tight.

\item Given any continuous function $f:\mbbr^{q+1}\to\mbbr^{q'}$ (for some $q'\in\mbbn$), we consider the function $H^{n}=(H^{1,n},H^{2,n})$ with
\begin{align}
H^{1,n}_{t}&:= 
\left(-Y_{[nt]/n},\,X_{[nt]/n}\right),
\nonumber \\
H^{2,n}_{t}&:=\frac{1}{n}\sum_{j=1}^{[nt]}f(Y_{t_{j-1}},X_{t_{j-1}}).
\nonumber
\end{align}
Then, we have 
$H^{1,n}\cip H^{1,\infty}:=\v2{(-Y,X)}$ in $\mcd(\mbbr^{q+1})$ 
and $H^{2,n}\cip H^{2,\infty}:=\int_{0}^{\cdot}f(Y_s,X_{s})ds$ in $\mcd(\mbbr^{q'})$, with which \eqref{LAD.AMN_add1} concludes the joint weak convergence in $\mcd(\mbbr^{2+q+q'})$:
\begin{equation}
(S^{n},H^{1,n},H^{2,n})\cil(w',H^{1,\infty},H^{2,\infty}).
\nn
\end{equation}
\end{itemize}
With these observations, we can apply \cite[VI.6.22]{JacShi03} to derive the weak convergence of stochastic integrals:
\begin{equation}
( H^{1,n}_{-}\cdot S^{n},H^{2,n} ) \cil ( H^{1,\infty}_{-}\cdot w',H^{2,\infty} ).
\nonumber
\end{equation}
which entails that, for any continuous function $f$,
\begin{align}
\left(\D_{0,n},\ \frac{1}{n}\sum_{j=1}^{n}f(Y_{t_{j-1}},X_{t_{j-1}})\right)
&\cil
\left(\frac{\sig_0^{-1}}{T} \int_{0}^{T}(-Y_{s},X_s) dw'_{s},\ \frac1T \int_{0}^{T} f(Y_s,X_{s})ds\right)
\nonumber \\
&\stackrel{\mcl}{=} \left(
\left\{\frac{\sig_0^{-2}}{T} \int_{0}^{T}
\begin{pmatrix}
Y_{t}^2 & -Y_{t}X_{t}^\top \\ -Y_{t}X_{t} & X_{t}^{\otimes 2}
\end{pmatrix}dt\right\}^{1/2} Z, ~\frac1T \int_{0}^{T} f(Y_s,X_{s})ds\right),
\nonumber
\end{align}
where $Z\sim N(0,1)$ independent of $\mcf$.
Now, by taking
\begin{equation}
f(x,y) = 
\frac{2\phi_{\beta_0}(0)}{\sig_{0}^2}
\begin{pmatrix}
y^2 & - x^\top y \\ - xy & x^{\otimes 2}
\end{pmatrix},
\nonumber
\end{equation}
we arrive at
\begin{equation}
\left(\D_{0,n}, \Gam_0\right)
\cil
\left(
\left\{\sig_0^{-2}\frac1T \int_{0}^{T}
\begin{pmatrix}
Y_{t}^2 & -Y_{t}X_{t}^\top \\ -Y_{t}X_{t} & X_{t}^{\otimes 2}
\end{pmatrix}dt\right\}^{1/2} Z, 
~\frac{2\phi_{\beta_0}(0)}{\sig_{0}^2}\frac1T \int_{0}^{T} 
\begin{pmatrix}
Y_{t}^2 & -Y_{t}X_{t}^\top \\ -Y_{t}X_{t} & X_{t}^{\otimes 2}
\end{pmatrix}
dt\right).
\nonumber
\end{equation}
In sum, applying Slutsky's theorem concludes that
\begin{equation}
\hat{w}_n = \Gam_{0}^{-1}\D_{0,n} + o_{p}(1) \cil 
MN_{q+1,\tz}\left(0,\, 
\frac{\sig_{0}^2}{4\phi_{\beta_0}(0)^2}
\left\{\frac1T \int_{0}^{T}
\begin{pmatrix}
Y_{t}^2 & -Y_{t}X_{t}^\top \\ -Y_{t}X_{t} & X_{t}^{\otimes 2}
\end{pmatrix}dt\right\}^{-1}\right).
\nonumber
\end{equation}

\bigskip

\noindent
{\bf Acknowledgement.}
The author should like to thank the anonymous reviewers for their detailed comments, based on which he could fix some essential mistakes and drastically improve the quality of the paper.
This work was partly supported by JSPS KAKENHI Grant Number 22H01139 and JST CREST Grant Number JPMJCR2115, Japan.

\medskip

\noindent
{\bf Conflict of interest.} 
The author declares that there is no conflict of interest.

\bigskip 
%

\begin{thebibliography}{10}

\bibitem{BorSch17}
S.~Borovkova and M.~D. Schmeck.
\newblock Electricity price modeling with stochastic time change.
\newblock {\em Energy Economics}, 63:51 -- 65, 2017.

\bibitem{BroMas18}
A.~Brouste and H.~Masuda.
\newblock Efficient estimation of stable {L}\'{e}vy process with symmetric
  jumps.
\newblock {\em Stat. Inference Stoch. Process.}, 21(2):289--307, 2018.

\bibitem{ChaSch12}
C.~J. Challis and S.~C. Schmidler.
\newblock A stochastic evolutionary model for protein structure alignment and
  phylogeny.
\newblock {\em Molecular biology and evolution}, 29(11):3575--3587, 2012.

\bibitem{CheKawMae03}
P.~Cheridito, H.~Kawaguchi, and M.~Maejima.
\newblock Fractional {O}rnstein-{U}hlenbeck processes.
\newblock {\em Electron. J. Probab.}, 8:no. 3, 14 pp. (electronic), 2003.

\bibitem{CleGlo20}
E.~Cl\'{e}ment and A.~Gloter.
\newblock Joint estimation for {SDE} driven by locally stable {L}\'{e}vy
  processes.
\newblock {\em Electron. J. Stat.}, 14(2):2922--2956, 2020.

\bibitem{Doo42}
J.~L. Doob.
\newblock The {B}rownian movement and stochastic equations.
\newblock {\em Ann. of Math. (2)}, 43:351--369, 1942.

\bibitem{DuM73}
W.~H. DuMouchel.
\newblock On the asymptotic normality of the maximum-likelihood estimate when
  sampling from a stable distribution.
\newblock {\em Ann. Statist.}, 1:948--957, 1973.

\bibitem{HjoPol93}
N.~L. Hj{\o}rt and D.~Pollard.
\newblock Asymptotics for minimisers of convex processes.
\newblock {\em Statistical Research Report, University of Oslo}, 1993.
\newblock Available at arxiv preprint arXiv:1107.3806, 2011.

\bibitem{HuLon09}
Y.~Hu and H.~Long.
\newblock Least squares estimator for {O}rnstein-{U}hlenbeck processes driven
  by {$\alpha$}-stable motions.
\newblock {\em Stochastic Process. Appl.}, 119(8):2465--2480, 2009.

\bibitem{Jac07}
J.~Jacod.
\newblock Asymptotic properties of power variations of {L}\'evy processes.
\newblock {\em ESAIM Probab. Stat.}, 11:173--196, 2007.

\bibitem{JacShi03}
J.~Jacod and A.~N. Shiryaev.
\newblock {\em Limit theorems for stochastic processes}, volume 288 of {\em
  Grundlehren der Mathematischen Wissenschaften [Fundamental Principles of
  Mathematical Sciences]}.
\newblock Springer-Verlag, Berlin, second edition, 2003.

\bibitem{Jeg82}
P.~Jeganathan.
\newblock On the asymptotic theory of estimation when the limit of the
  log-likelihood ratios is mixed normal.
\newblock {\em Sankhy\=a Ser. A}, 44(2):173--212, 1982.

\bibitem{Jeg95}
P.~Jeganathan.
\newblock Some aspects of asymptotic theory with applications to time series
  models.
\newblock {\em Econometric Theory}, 11(5):818--887, 1995.
\newblock Trending multiple time series (New Haven, CT, 1993).

\bibitem{JhwMar14}
D.-C. Jhwueng and V.~Maroulas.
\newblock Phylogenetic {O}rnstein-{U}hlenbeck regression curves.
\newblock {\em Statist. Probab. Lett.}, 89:110--117, 2014.

\bibitem{Leh99}
E.~L. Lehmann.
\newblock {\em Elements of large-sample theory}.
\newblock Springer Texts in Statistics. Springer-Verlag, New York, 1999.

\bibitem{Mas10ejs}
H.~Masuda.
\newblock Approximate self-weighted {LAD} estimation of discretely observed
  ergodic {O}rnstein-{U}hlenbeck processes.
\newblock {\em Electron. J. Stat.}, 4:525--565, 2010.

\bibitem{Mas19spa}
H.~Masuda.
\newblock Non-{G}aussian quasi-likelihood estimation of {SDE} driven by locally
  stable {L}\'{e}vy process.
\newblock {\em Stochastic Process. Appl.}, 129(3):1013--1059, 2019.

\bibitem{PKAS11}
M.~Perninge, V.~Knazkins, M.~Amelin, and L.~S{\"o}der.
\newblock Modeling the electric power consumption in a multi-area system.
\newblock {\em European Transactions on Electrical Power}, 21(1):413--423,
  2011.

\bibitem{SamKni09}
D.~M.~M. Samarakoon and K.~Knight.
\newblock A note on unit root tests with infinite variance noise.
\newblock {\em Econometric Rev.}, 28(4):314--334, 2009.

\bibitem{Swe80}
T.~J. Sweeting.
\newblock Uniform asymptotic normality of the maximum likelihood estimator.
\newblock {\em Ann. Statist.}, 8(6):1375--1381, 1980.
\newblock Corrections: (1982) {\it Annals of Statistics} {\bf 10}, 320.

\bibitem{VAKB19}
H.~Verdejo, A.~Awerkin, W.~Kliemann, and C.~Becker.
\newblock Modelling uncertainties in electrical power systems with stochastic
  differential equations.
\newblock {\em International Journal of Electrical Power \& Energy Systems},
  113:322 -- 332, 2019.

\bibitem{Zac71}
S.~Zacks.
\newblock {\em The theory of statistical inference}.
\newblock John Wiley \& Sons, Inc., New York-London-Sydney, 1971.
\newblock Wiley Series in Probability and Mathematical Statistics.

\bibitem{ZhaZha13}
S.~Zhang and X.~Zhang.
\newblock A least squares estimator for discretely observed
  {O}rnstein-{U}hlenbeck processes driven by symmetric {$\alpha$}-stable
  motions.
\newblock {\em Ann. Inst. Statist. Math.}, 65(1):89--103, 2013.

\end{thebibliography}

\def\cprime{$'$} \def\polhk#1{\setbox0=\hbox{#1}{\ooalign{\hidewidth
  \lower1.5ex\hbox{`}\hidewidth\crcr\unhbox0}}} \def\cprime{$'$}
  \def\cprime{$'$}

\end{document}